\documentclass{amsart}

\usepackage{fullpage}

\numberwithin{equation}{section}

\usepackage{amssymb}
\usepackage{enumerate}
\usepackage{changebar}
\newtheorem{thm}{Theorem}[section]
\newtheorem{lem}[thm]{Lemma}
\newtheorem{cor}[thm]{Corollary} 
\newtheorem{sublem}[thm]{Sub-lemma}
\newtheorem{prop}[thm]{Proposition}

\newtheorem{exmp}{Example}
\newtheorem{rem}[thm]{Remark}

\newcommand\Co{{\mathcal C}}

\newcommand\F{{\mathcal F}}
\newcommand\Ka{{\mathcal K}}
\newcommand\Lp{{\mathcal L}}
\newcommand\Or{{\mathcal O}}
\newcommand\M{{\mathcal M}}

\newcommand\N{{\mathbb N}}
\newcommand\R{{\mathbb R}}
\newcommand\To{{\mathbb T}}
\newcommand\Z{{\mathbb Z}}

\newcommand\ve{\varepsilon}

\newcommand\vf{\varphi}

\newcommand\Id{\text{\bf Id}}

\newcommand{\cs}{\operatorname{C}}

\begin{document}

\title[Convergence to equilibrium for intermittent symplectic
maps]{Convergence to equilibrium for intermittent symplectic maps}  

\author{Carlangelo Liverani, Marco Martens}
\address{Carlangelo Liverani\\
Dipartimento di Matematica\\
II Universit{\`a} di Roma (Tor Vergata)\\
Via della Ricerca Scientifica, 00133 Roma, Italy.}
\email{{\tt liverani@mat.uniroma2.it}} 
\address{Marco Martens\\
University of Groningen, Department of Mathematics
P.O. Box 800, 9700 AV Groningen, The Netherlands}
\email{{\tt m.martens@math.rug.nl}} 
\date{\today}
\thanks{One of us (C.L.) would like to thank S.Vaienti and M.Benedicks
for helpful discussions. In addition, we acknowledge the support from
the ESF Program  PRODYN, I.B.M., M.I.U.R. and NSF} 
\begin{abstract}
We investigate a class of area preserving non-uniformly hyperbolic
maps of the two torus. First we establish some results on the
regularity of the invariant foliations, then we use this knowledge
to estimate the rate of mixing.
\end{abstract}
\maketitle

%VT
\thispagestyle{empty} 
\def\IMSmarkvadjust{0 pt}
\def\IMSmarkhadjust{0 pt}
\def\IMSmarkhpadding{0 pt}
\def\IMSpubltext{Published in modified form:}
\def\SBIMSMark#1#2#3{
 \font\SBF=cmss10 at 10 true pt
 \font\SBI=cmssi10 at 10 true pt
 \setbox0=\hbox{\SBF \hbox to \IMSmarkhpadding{\relax}
                Stony Brook IMS Preprint \##1}
 \setbox2=\hbox to \wd0{\hfil \SBI #2}
 \setbox4=\hbox to \wd0{\hfil \SBI #3}
 \setbox6=\hbox to \wd0{\hss
             \vbox{\hsize=\wd0 \parskip=0pt \baselineskip=10 true pt
                   \copy0 \break%
                   \copy2 \break% 
                   \copy4 \break}}
 \dimen0=\ht6   \advance\dimen0 by \vsize \advance\dimen0 by 8 true pt
                \advance\dimen0 by -\pagetotal
	        \advance\dimen0 by \IMSmarkvadjust
 \dimen2=\hsize \advance\dimen2 by .25 true in
	        \advance\dimen2 by \IMSmarkhadjust

%
%   Check for publication info
%
%  \newread\jref
  \openin2=publishd.tex
  \ifeof2\setbox0=\hbox to 0pt{}
  \else 
     \setbox0=\hbox to 3.1 true in{
                \vbox to \ht6{\hsize=3 true in \parskip=0pt  \noindent  
                {\SBI \IMSpubltext}\hfil\break
                \input publishd.tex 
                \vfill}}
  \fi
  \closein2
  \ht0=0pt \dp0=0pt
 \ht6=0pt \dp6=0pt
 \setbox8=\vbox to \dimen0{\vfill \hbox to \dimen2{\copy0 \hss \copy6}}
 \ht8=0pt \dp8=0pt \wd8=0pt
 \copy8
 \message{*** Stony Brook IMS Preprint #1, #2. #3 ***}
}

\SBIMSMark{2005/01}{January 2005}{}

\section{Introduction}
\label{sec:intro}
In recent years the physics community has devoted an increasing
attention to anomalous properties of physical systems (e.g., anomalous
transport, anomalous diffusion, anomalous conductivity, etc.). Such
properties have proven 
relevant in many fields such as thermal conductivity, kinetic
equations, plasma physics, etc. and they are widely believed to be
dynamical in nature. In fact, such phenomena seem to depend on the
weak chaotic properties of the underling dynamics, see \cite{Za} and
references therein for a detailed discussion. The 
basic idea is that, while uniformly hyperbolic dynamics gives rise to
normal transport properties (consider for example the diffusive
behavior in a finite horizon Lorentz gas \cite{BuSi})
non-uniform hyperbolicity gives rise to different behavior
(e.g. the anomalous diffusion believed to occur in infinite horizon
Lorentz gas \cite{Bl}) due to weaker
mixing properties (e.g. polynomial decay of correlations) and regions
in which the motions is rather regular and where the systems spend a
substantial fraction of time ({\em sticky} regions).  

Unfortunately, the theoretical understanding of dynamical models with
polynomial decay of correlations is extremely limited, hence the
necessity to rigorously investigate relevant toy models. The only well
understood cases are expanding one dimensional maps with a neutral
fixed point. Such maps were proposed as a model of intermittent
behavior in fluids (\cite{PM}) and have been widely studied. It has been
proven that such maps enjoy polynomial decay of correlations with the
rate depending on the behavior of the fixed point \cite{LSV, Yo1,
Hu1, Hu2, Sa, Go1}. In addition, when the
decay of correlation is sufficiently slow, the observables do not
satisfy the Central Limit Theorem or the Invariance Principle but
rather, when properly rescaled, some stable law (\cite{Zw, Go2, PoSh}).

Some, more partial, results exist for multidimensional expanding
maps \cite{PoYu} as well. Yet, the usual physically relevant models are
connected to Hamiltonian dynamics and, to our knowledge, no rigorous
results are available in such a situation. The simplest case which
retain some Hamiltonian flavor is clearly a two dimensional area
preserving map. In fact, mixing area preserving maps of the two
dimensional torus with a neutral fixed point (the simplest type of
sticky set) have been investigated
numerically \cite{ArPr} predicting the possibility of a polynomial decay
of correlations. 
In this paper we consider a class of non-uniformly hyperbolic symplectic maps \eqref{eq:map}
of the two torus where the hyperbolicity breaks down because of a non-hyperbolic fixed point. In fact, the linearized
dynamics at the fixed point is a shear \eqref{eq:der}. We prove that
the decay of correlations is polynomial and, more precisely, decays at
least as $n^{-2}$ and, in some sense, one cannot expect much
more. Note that the example treated numerically in \cite{ArPr} is a
special case of the present setting. In \cite{ArPr} the predicted
decay was $n^{-2.5}$. This emphasizes the difficulty to investigate
such issues and the strong need for more theoretical work on the subject.

The result in the present paper is based on a precise quantitative
analysis of the angle between the stable and the unstable
direction. This angle turns out to degenerate approaching the origin
(where the non hyperbolic fixed point is located). Once such a control
is achieved it is possible to obtain a bound on the expansion and
contraction in the system. Such expansion turns out to be only
polynomial, in contrast with the uniformly hyperbolic case where it is
exponential. In turn, the bound on the expansion allows to study the
regularity of the stable and unstable foliation. It turns out that
they are $\Co^1$ away from the origin. This suffice to apply a simple
random approximation technique that allows estimating the speed of the
correlations. 

As the rate of convergence to equilibrium is of order $n^{-2}$, see Theorem
\ref{thm:main}, the
Central Limit Theorem holds for zero average observable, see Corollary 
\ref{cor:clt},
so the model does not exhibit anomalous statistical behavior in this
respect. Yet, it clearly exhibits an intermittent behavior and it shows the
mechanism whereby slow decay of correlations may arise. The present
work emphasizes the need to carry out similar studies in cases where the
set producing intermittency has a more complex structure than a simple
isolated point.

The paper is organized as follows: section \ref{sec:results} details
the model and makes precise the results. Section \ref{sec:fix} studies
the local dynamics at the fixed point and, in particular the
properties of its stable and unstable manifolds. This can be achieved
in many way, here we find most efficient to apply a variational
technique. Section \ref{sec:narrow} establishes a precise bound for
the angle between the stable and the unstable direction at each
point. As anticipated, such a bound yields an a priory bound on the
expansion and contractions rates in the systems, these are obtained in
section \ref{sec:expansion}. The latter result suffices to apply
standard distortion estimates that, in turn, allow to prove precise
results on the regularity of the invariant foliation and the
holonomies, see section \ref{sec:regularity} and section
\ref{sec:holo} respectively.  Next, in section
\ref{sec:random} we introduce a random perturbation of the above map
and investigate its statistical properties that, thanks to the added
randomness, can be addressed fairly easily. The relevance of the above
random perturbation is that the limit of zero noise allows to easily obtain a
bound on the rate of mixing in the original map, we do this in section
\ref{sec:decay}. Finally, in section \ref{sec:lower}, we show that the
obtained bound is close to being optimal. The paper ends with Remark
\ref{rem:problems} pointing to the unsatisfactory nature of some of
the present results and the need to investigate the related open
problems.

\section{The model and the results}
\label{sec:results}
For each $h\in\Co^{\infty}(\To^1,\To^1)$ we define the map
$T:\To^2\to\To^2$ by\footnote{\label{foot:uno} Note that the following formula is
equivalent, by the symplectic change of variable $q=x-y$, $p=y$, to
the map
\[
\widetilde T(q,p)=\begin{cases}
              q+p&\quad\mod 1\\
              p+h(q+p)&\quad \mod 1
              \end{cases}
\]
which belongs to the standard map family. Yet, the functions $h$
considered here differ substantially from the sine function which would
correspond to the classical Chirikov-Taylor well known example.}
\begin{equation}\label{eq:map}
T(x,y)=\begin{cases}
              x+h(x)+y&\quad\mod 1\\
              h(x)+y&\quad\mod 1
              \end{cases}
\end{equation} 
We moreover require the following properties
\begin{enumerate}
\item $h(0)=0$  (zero is a fixed point);
\item $h'(0)=0$ (zero is a neutral fixed point)
\item $h'(x)>0$ for each $x\neq 0$ (hyperbolicity)
\end{enumerate}
Note that conditions (2--3) imply that zero is a minimum for $h'$, which forces 
\[
h''(0)=0;\quad h'''(0)\geq 0.
\]
We will restrict to the generic case
\begin{enumerate}
\item[(4)] $h'''(0)>0$.
\end{enumerate}
In order to simplify the discussion we will also assume the following
symmetry
\begin{enumerate}
\item[(5)] $h(-x)=-h(x)$.
\end{enumerate}

This means that we can write
\begin{equation}
\label{eq:h}
h(x)=bx^3+\Or(x^5).
\end{equation}
\begin{rem}
Note that two facts implied by the above assumptions are not necessary and could be
done away with at the price of more extra work: the hypothesis that
there is only one neutral fixed point (finitely many neutral periodic
orbits would make little difference) and the symmetry (5). We assume
such facts only to simplify the presentation of the arguments.
\end{rem}
Since the derivative of the map is given by
\begin{equation}\label{eq:der}
DT=\begin{pmatrix}
   1+h'(x)&1\\
   h'(x)&1
   \end{pmatrix}
\end{equation}
$\det(DT)=1$, thus the Lebesgue measure $m$ is an invariant measure (the
maps are symplectic). From now on we will consider the dynamical
systems $T: (\To^2,m)\rightarrow (\To^2,m)$. 

Formula \eqref{eq:der} and property (3) imply that 
the cone $\Co_+=\{v\in\R^2\;|\;Q(v):=\langle v_1,\,v_2\rangle\geq 0\}$
is invariant for $DT$. In additions, it is easy to check that $D_\xi
T^2 \Co_+\subset \hbox{int }\Co_+\cup\{0\}$ for all
$\xi\in\To^2\backslash \{0\}$. From
this and the general theory, see \cite{LW}, follows immediately

\begin{thm}
The above described dynamical systems are non-uniformly hyperbolic and mixing.
\end{thm} 
\begin{exmp}
An interesting concrete example for the above setting is given by the function 
$h(x):=x-\sin x$.
\end{exmp}
The question remains about the rate of mixing, this is the present
topic. 
\begin{rem}
In the following by $\cs$ we designate a generic
constant depending only on $T$. Accordingly, its value may vary from
an occurrence to the next. In the instances when we will need a constant
of the above type but with a fixed value we will use sub-superscripts.
\end{rem}
\begin{thm}\label{thm:main}
For each $f,g\in\Co^{1}(\To^2,\R)$, $\int f=0$, holds
true\footnote{In fact, a slightly sharper bound holds, see \eqref{eq:sharp}.}
\[
\left|\int f g\circ T^n\right|\leq \cs
\|f\|_{\Co^{1}}\|g\|_{\Co^{1}}n^{-2}(\ln n)^4.
\]
\end{thm}
\begin{rem}
As in other similar cases \cite{LSV, L, Po} the logarithmic
correction is almost certainly 
an artifact of the technique of the proof. It could probably be
removed by using a more sophisticated (and thus more technically
involved) approach. See also section \ref{sec:lower}.
\end{rem}
Form Theorem \ref{thm:main} many facts follow, just to give an
example let us mention the following result that can be obtained from
Theorem 1.2 in \cite{L3}.
\begin{cor}[CLT]\label{cor:clt}
Given $f\in\Co^1$, $\int f=0$, the random variable
\[
\frac 1{\sqrt n}\sum_{i=0}^{n-1}f\circ T^i
\]
converges in distribution to a Gaussian variable with zero mean and
finite variance $\sigma$. In addition, $\sigma=0$ iff there exists
$\vf\in L^1$ such that $f=\vf-\vf\circ T$.\footnote{In particular this means
that the average of $f$ on each periodic orbit must be zero.}
\end{cor}
The rest of the paper is devoted to the proof of
Theorem \ref{thm:main} that will find its conclusion in section
\ref{sec:decay}. The basic 
fact needed in the proof, a fact of independent 
interest and made quantitatively precise in Lemma \ref{lem:distreg1}, is
the following. 
\begin{thm}\label{thm:mainlem}
The stable and unstable distributions are $\Co^1$ in $\To^2\backslash\{0\}$.
\end{thm}

\section{The fixed point manifolds}\label{sec:fix}

As usual we start by studying the local dynamics near the fixed
point. The first basic fact is the existence of stable and unstable
manifolds. This is rather standard, yet since we need some
quantitative information we will construct them explicitly.

Instead of constructing them via usual fixed point arguments it turns
out to be faster to use a variational method.

\subsection{A variational argument}

Let us consider, in a neighborhood of zero, the function
\begin{equation}\label{eq:generating}
L(x,x_1):=\frac 12 (x-x_1)^2+G(x)\;;\quad G(x):=\int_0^xh(z) dz.
\end{equation}
By setting
\[
\begin{split}
&y:=-\frac{\partial L}{\partial x}=x_1-x-h(x)\\
&y_1:=\frac{\partial L}{\partial x_1}=x_1-x
\end{split}
\]
we have $(x_1,y_1)=T(x,y)$, that is {\sl $L$ is a generating function
for the map \eqref{eq:map}}.

Then, for each $a\in\R$, we define the Lagrangian $\Lp_a:\ell^2(\N)\to\R$ by
\begin{equation}\label{eq:lagrange}
\Lp_a(x):=\sum_{n=1}^\infty L(x_n,x_{n+1})+L(a,x_1).
\end{equation}
The justification of the above definition rests in the following
Lemma.
\begin{lem}\label{lem:critical}
For each $a\in(-1,1)$, holds true $\Lp_a\in\Co^1(\ell^2(\N))$. In addition, if
$x\in \ell^2(\N)$ is such that $D_x\Lp_a=0$, then setting $x_0=a$ and 
$y_n:=x_{n+1}-x_n-h(x_n)$, we have $T^n(x_0,y_0)=(x_n,y_n)$.
\end{lem}
\begin{proof}
First of all \eqref{eq:h} implies that there exists $\cs>0$ such that
$|G(x)|\leq \cs x^4$. It is then easy to see that $\Lp_a$ is well defined
for each sequence in $\ell^2(\N)$.

Next, for each $n\in\N$ let us define
$(\nabla\Lp_a)_n:=\partial_{x_n}\Lp_a$.
Clearly,
\[
\begin{split}
(\nabla\Lp_a)_1(x)&=2x_1-x_2+h(x_1)-a\\
(\nabla\Lp_a)_n(x)&=2x_n-x_{n+1}-x_{n-1}+h(x_n).
\end{split}
\]
Of course, for $x\in\ell^2(\N)$, $ \nabla\Lp_a(x)\in\ell^2(\N)$, it is
then trivial to check 
that $D_x\Lp_a(v)=\langle \nabla\Lp_a(x), v\rangle$.
The last statement follows by a direct computation.
\end{proof}

By the above Lemma it is clear that one can obtain the stable
manifolds of the fixed point from the critical points of $\Lp_a$, it
remains to prove that such critical points do exist. We will start by
considering the case $a\geq 0$.

Define
\begin{equation}\label{eq:fixedset}
Q_B:=\{x\in\ell^2(\N)\;|\;\,|x_n-\frac A{n+c}|\leq B(n+c)^{-\frac
32}\}
\end{equation}
where $A:=\sqrt{\frac 2b}$; $c:=\frac Aa$.

It is immediate to check that $Q_B$ is compact and convex. In
addition, if $a$ is sufficiently small, then $G$ is strictly convex on
$[-2a,2a]$ which implies that $\Lp_a|_{Q_B}$ is strictly
convex. Accordingly, $\Lp_a$ has minimum in $Q_B$, moreover the strict
convexity implies that such a minimum is unique, for $a$ fixed.

Let us call $x(a)$  the point in $Q_B$ where $\Lp_a$ attains its
minimum.
\begin{lem}\label{lem:stab}
For $a$ small enough, $D_{x(a)}\Lp_a=0$.
\end{lem}
\begin{proof}
Suppose that $\partial_{x_n}\Lp_a(x(a))\neq 0$ for some $n\in\N$, for
example suppose it is negative. Then $x(a)$ is on the border of $Q_B$,
say $x(a)_n=A(n+c)^{-1}+B(n+c)^{-\frac 32}$, otherwise we could increase
$x(a)_n$ and decrease $\Lp_a$ still remaining in $Q_B$, contrary to
the assumption. But then
\[
\begin{split}
\partial_{x_n}\Lp_a(x(a))&=2x(a)_n-x(a)_{n-1}-x(a)_{n+1}+h(x(a)_{n})\\
&=\frac {2A}{n+c}+\frac{2B}{(n+c)^{\frac
32}}-x(a)_{n-1}-x(a)_{n+1}+h(\frac {A}{n+c}+\frac{B}{(n+c)^{\frac
32}}) \\
&\geq -\frac{2A}{[(n+c)^2-1](n+c)}+\frac
{2A}{(n+c)^3}+\frac{3bA^2B}{(n+c)^\frac 72}+\Or((n+c)^{-4})\\
&+\frac{B(n+c)^3\{2(1-(n+c)^{-2})^{\frac 32}-(1-(n+c)^{-1})^{\frac
32} -(1+(n+c)^{-1})^{\frac 32}}{(n+c)^{\frac 32}[(n+c)^2-1]^{\frac
32}}\\
&=\frac{(6-\frac {15}4)B}{(n+c)^{\frac 72}}+\Or((n+c)^{-4})\geq 0
\end{split}
\]
provided $a$ is sufficiently small. We have thus a contradiction. The
other possibilities are analyzed similarly.
\end{proof}

To conclude we need some information on the regularity of $x(a)$ as a
function of $a$. Unfortunately, the
implicit function theorem does not applies since $D^2\Lp_a$ does not
have a spectral gap, yet for our purposes a simple estimate
suffices.
\begin{lem}\label{lem:lipman}
$x_1(a)$ is a Lipschitz function of $a$. Moreover, when derivable,
\[
|y_0(a)'|\leq \cs |x(a)_0|.
\]
\end{lem}
\begin{proof}
By Lemma \ref{lem:stab} it follows, for each $a,a'$ sufficiently small
\[
\partial_{x_n}\Lp_{a'}(x(a'))=\partial_{x_n}\Lp_{a}(x(a))=0
\]
that is
\[
\partial_{x_n}\Lp_{a'}(x(a'))-\partial_{x_n}\Lp_{a'}(x(a))=\partial_{x_n}\Lp_{a}(x(a))
-\partial_{x_n}\Lp_{a'}(x(a))
\]
which yields
\begin{equation}\label{eq:strange}
\begin{split}
&(2+h'(\xi_1))\zeta_1-\zeta_2=a'-a\\
&(2+h'(\xi_n))\zeta_n-\zeta_{n+1}-\zeta_{n-1}=0,
\end{split}
\end{equation}
where $\zeta_n:=x(a')_n-x(a)_n$ and $\xi_n\in[x(a)_n,x(a')_n]$.

Notice that, if $|\zeta_n|\geq|\zeta_{n-1}|$, then
\[
|\zeta_{n+1}|=|(2+h'(\xi_n))\zeta_n-\zeta_{n-1}|\geq
2|\zeta_n|-|\zeta_{n-1}|\geq |\zeta_n|.
\]

Thus, by induction, if $|\zeta_n|\geq |\zeta_{n-1}|$, then $|\zeta_m|\geq
|\zeta_{n-1}|$ for each $m\geq n$, which would imply $\zeta_{n-1}=0$ since
$\zeta\in\ell^2(\N)$. But then
$(2+h'(\xi_n))\zeta_n=\zeta_{n+1}$, that is
$|\zeta_{n+1}|\geq|\zeta_n|$. Accordingly, again by induction,
$\zeta_m=0$ for each $m\geq n-1$. This means that we can restrict
ourselves to the case $\zeta_n\neq 0$, $|\zeta_n|\geq |\zeta_{n+1}|$.
Hence,
\[
|a'-a|=|(2+h'(\xi_1))\zeta_1-\zeta_2|\geq 2|\zeta_1|-|\zeta_2|\geq
|\zeta_1|.
\]
That is $|x(a')_n-x(a)_n|\leq |x_1(a')-x_1(a)|\leq |a'-a|$.

Finally, summing \eqref{eq:strange} over $n\in\N$, 
\[
\sum_{n=1}^\infty h'(\xi_n)\zeta_n=-\zeta_1+a'-a.
\]
Thus, where all the $x(a)_n$ are differentiable (a full measure set),
$|x(a)_n'|\leq|x(a)_1'|\leq 1$ and
\begin{equation}\label{eq:derivxy}
\begin{split}
x(a)_1'&=1-\sum_{n=1}^\infty h'(x(a)_n)x(a)_n'\\
y(a)_0'&=-\sum_{n=0}^\infty h'(x_n)x(a)_n'.
\end{split}
\end{equation}
Accordingly,
\[
|y_0(a)'|\leq 6b\sum_{n=0}^\infty x_n^2\leq \cs a.
\]
\end{proof}

Clearly, the above Lemma implies that, calling $(x,\gamma_s(x))$ the
graph of the stable manifold, $\gamma_s\in Lip(-1,1)$. The case $a\leq
0$ and the unstable manifolds can be treated similarly, yet there
exists a faster--and more
instructive--way.

\subsection{Reversibility}

Notice that the map $T$ is {\sl reversible}
with respect to the transformations\footnote{While the reversibility for
$\Pi$ is a general fact, the one for $\Pi_1$ depends on the simplifying
symmetry hypothesis (5).}
\begin{equation}\label{eq:invol}
\Pi(x,y):=(x, -y-h(x)); \quad \Pi_1(x,y):=(-x, y+h(x))
\end{equation}
Indeed, $\Pi^2=\Pi_1^2=\Id$ and $\Pi T\Pi=\Pi_1 T\Pi_1=T^{-1}$. 

\begin{rem}\label{rem:rev}
The reversibility implies that, for $x\geq 0$,
$(x,\gamma_u(x))=\Pi(x,\gamma_s(x))$, and, for $x\leq 0$,
$(x,\gamma_u(x))=\Pi_1(-x,\gamma_s(-x))$ is the unstable manifold of
zero. 
\end{rem}

\subsection{A quasi-Hamiltonian}

\begin{figure}[ht]
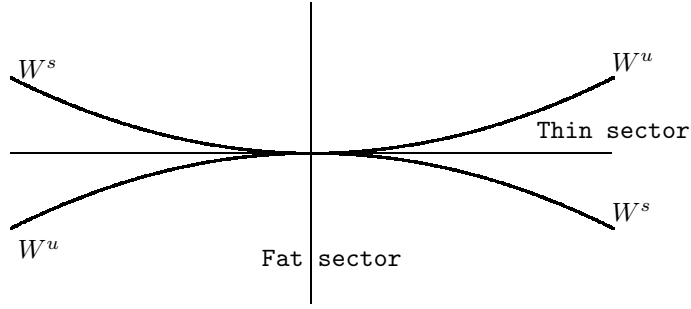
\
\centering
\put(-40,0){\line(1,0){80}}
\put(0,-20){\line(0,1){40}}
\put(-6.5,-15){{\tt Fat sector}}
\put(30,2){{\tt Thin sector}}
\thicklines
\qbezier(-40, 10)(0,-10)(40,10)
\qbezier(-40, -10)(0,10)(40,-10)
\put(-39,10){$W^s$}
\put(-39,-14){$W^u$}
\put(40,11){$W^u$}
\put(40,-9){$W^s$}
\put(0,-24){\ }
\caption{The manifolds of the fixed point}\label{fig:man}
\end{figure}

To study the motion near the fixed point it is helpful to find a local
``Hamiltonian'' function. By Hamiltonian function we mean a function
that is locally invariant for the dynamics. Such a function can be
computed as a formal power series starting by the relation $H\circ
T=H$. In fact, we are interested only in a suitable approximation. A
direct computation yields that, by defining $G(x):=\int_0^x h(z)dz$ and
\begin{equation}\label{eq:hamiltonian}
H(x,y):=\frac 12 y^2-G(x)+\frac 12 h(x)y-\frac 1{12}h'(x)y^2+\frac 1 {12}h(x)^2,
\end{equation}
holds true\footnote{In fact, setting $(x_1,y_1):=T(x,y)$, holds
\[
\begin{split}
&H(x_1,y_1)-H(x,y)=yh(x)+\frac 12 h(x)^2-h(x)y-h(x)^2-\frac
12h'(x)y_1^2-\frac 16 h''(x)y_1^3+\frac 12h'(x)y_1^2\\
&+\frac 12 h(x)^2
+\frac 14 h''(x)y_1^3-\frac 1{12}h''(x)y_1^3-\frac
16h'(x)h(x)y_1+\frac 16 h'(x)h(x)y_1+\Or(x^8+x^4y^2+y^4).
\end{split}
\]
}
\begin{equation}\label{eq:appham}
H(T(x,y))-H(x,y)=\Or (x^8+y^4).
\end{equation}
This approximate conservation law suffices to obtain rather precise
information on the near fixed point dynamics.\footnote{The reader
  should be aware that it is possible to do much better, that is to
  obtain an exponentially precise conservation law, see \cite{L2},
  \cite{BG}.} 
The first application 
is given by the following information on the stable manifold.
\begin{lem}\label{lem:manifzero}
For $x\geq 0$ sufficiently small holds
\[
\gamma_s(x)= -A^{-1}x^2+\Or(x^3).
\] 
\end{lem}
\begin{proof}
Using the notation of Lemma \ref{lem:stab}, for fixed $a$ we get
\begin{equation}\label{eq:unoy}
y_n=x_{n+1}-x_n-h(x_n)=\Or((n+c)^{-\frac 32}).
\end{equation}
Hence
\[
H(x_n,y_n)=\Or((n+c)^{-3}).
\]
Using equation \eqref{eq:appham} we have
\[
H(x_0,y_0)=H(x_n,y_n)+\Or(\sum_{i=0}^{n-1}x_i^8+y_i^4)
\]
that is
\[
|H(x_0,y_0)|\leq \cs \left\{(n+c)^{-3}+
\sum_{i=0}^{\infty}(n+c)^{-6}\right\}\leq \cs
\left\{(n+c)^{-3}+c^{-5}\right\}.
\]
Taking the limit for $n$ to infinity in the above expression and
remembering the definition of $c$ follows
\[
|H(x_0,y_0)|\leq\cs x_0^{5}.
\]
Since equation \eqref{eq:unoy} implies $|y_0|=\Or(x_0^{\frac 32})$,
from \eqref{eq:hamiltonian} we have
\[
y_0^2-\frac b2 x_0^4+bx_0^3y_0=\Or(x_0^5)
\]
from which the lemma follows.
\end{proof}
According to Lemma \ref{lem:manifzero}, the local picture of the
manifolds is given by Figure \ref{fig:man}.

\subsection{Manifold regularity}
Since in the previous section we have seen that the manifold are
Lipschitz curves, we can define the dynamics restricted to the unstable
manifold:
\begin{equation}\label{eq:fu}
f_u(x):=x+h(x)+\gamma_u(x).
\end{equation}

Our next task is to obtain sharper information on the manifolds
regularity.\footnote{Of course, the manifold should be as smooth as $h$, but
this results is not needed in the following while we do need an
explicit bound on the curvature.}
\begin{lem}\label{lem:manifoldreg}
The unstable manifold of the fixed point is $\Co^2$, apart from zero.
\end{lem}
\begin{proof}
It is clearly enough to show that $\gamma_u\in \Co^2$ apart from
zero. To do so call $u(x)=\gamma_u'(x)$ (the derivative exists almost
everywhere since $\gamma_u$ is Lipschitz). The tangent vector to the
unstable manifold has the form $(1,u)$. On the other hand
\[
\begin{split}
\lambda_u(f_u^{-1}(x))
\begin{pmatrix}1\\u(x)\end{pmatrix}&:=\begin{pmatrix}
                                    1+h'(f_u^{-1}(x))&1\\
                                     h'(f_u^{-1}(x))&1
                                    \end{pmatrix}
                          \begin{pmatrix}1\\u(f_u^{-1}(x))\end{pmatrix}\\
&=:\lambda_u(f_u^{-1}(x))\begin{pmatrix}1\\F(f_u^{-1}(x),u(f_u^{-1}(x)))
\end{pmatrix}
\end{split}
\]
where
\begin{equation}\label{eq:unsdef}
\begin{split}
&\lambda_u(x):=1+h'(x)+u(x)=f'_u(x)\\
&F(x,u):=1-\frac{1}{1+h'(x)+u}.
\end{split}
\end{equation}
Accordingly, setting $x_i:=f_u^{-i}(x)$, holds
$u(x_i)=F(x_{i+1},u(x_{i+1}))$. Next, let $v_i:=2A^{-1}x_i$, then a
direct computation yields $F(x_i,v_i)-v_{i-1}=\Or(x_i^3)$, thus
\[
|u(x_i)-v_i|\leq |F(x_{i+1},u(x_{i+1}))-F(x_{i+1}, v_{i+1})|+\cs x_i^3
\leq |u(x_{i+1})- v_{i+1}|+\cs x_i^3.
\]
By induction, and equation \eqref{eq:fixedset}, it follows
$|u(x)-2A^{-1}x|\leq \cs\sum_{i=0}^\infty x_i^3\leq \cs x^2$.

On the other hand, given a different
point $z$, it holds
\[
\begin{split}
u(x)-u(z)&=F(x_1,u(x_1))-F(z_1,u(z_1))\\
&=\lambda_u(x_1)^{-1}\lambda_u(z_1)^{-1}
\left\{(u(x_1)-u(z_1))+h'(x_1)-h'(z_1))\right\}.
\end{split}
\]
Iterating the above equation yields
\begin{equation}\label{eq:unstreg}
\begin{split}
u(x)-u(z)=&\lambda_{u,n}(x)^{-1}\lambda_{u,n}(z)^{-1}(u(x_n)-u(z_n))\\
&\ \ +\sum_{k=1}^{n}\lambda_{u,k}(x)^{-1}\lambda_{u,k}(z)^{-1}(h'(x_k)-h'(z_k))
\end{split}
\end{equation}
where $\lambda_{u,n}(x):=\prod_{k=1}^{n}\lambda_u(f_u^{-k}(x))$.
Next, let $x=a$, then, accordingly to Lemma
\ref{lem:stab}, equation \eqref{eq:fixedset} and Remark \ref{rem:rev}, we have
$x_i= A(i+c)^{-1}+\Or((i+c)^{3/2})$. This means that, in a
sufficiently small neighborhood of zero, and for $z$ sufficiently
close to $x$ holds
\[
\begin{split}
\lambda_{u,n}(x)&\geq
\prod_{k=1}^{n}\left(1+v_k-\cs x_i^2\right)
\geq \prod_{k=1}^{n}\left(1+\frac 2{k+c}-\cs (k+c)^{-3/2}\right)\\
&\geq e^{\sum_{k=1}^{n}\frac{2}{k+c}-\cs (k+c)^{-3/2}}
\geq \cs (A^{-1}an+1)^{2}
\end{split}
\]
provided $x$ is close  
enough to zero. The same estimate holds for $\lambda_{u,n}(z)$. 

This implies that $u$ is continuous. Indeed, for each $\ve>0$ choose
$n_0(x)\in\N$ such that $\lambda_{u,n_0(x)}(x_n)^{-1}\leq \ve$, then
\[
|u(x)-u(z)|\leq \sum_{k=1}^{n_0(x)}|h'(x_k)-h'(z_k)|+\ve 
\]
and we can thus choose $z$ close enough to $x$ such that
$|u(x)-u(z)|\leq 2\ve$. Note that this implies the continuity of the
$\lambda_{u,n}$ as well.

To conclude, we choose $n(z)$ such that $\cs(A^{-1}an(z)+1)^{-4}\geq
|x-z|^{1+\alpha}$, for some $\alpha>0$, 
accordingly 
\[
u(x)-u(z)=\sum_{k=1}^{n(z)}\lambda_{u,k}(x)^{-1}
\lambda_{u,k}(z)^{-1}(h'(x_k)-h'(z_k))+{\mathcal
O}(|x-z|^{1+\alpha})
\]
Since the series is uniformly convergent we have
\begin{equation}\label{eq:unstder}
u'(x)=\sum_{k=1}^{\infty}\lambda_{u,k}(x)^{-3}h''(x_k)
\end{equation}
from which the lemma follows.\footnote{Remark that to obtain the 
result on a larger neighborhood it suffices to iterate the unstable  
manifold forward.}
\end{proof}

\begin{rem}
All the above results for the unstable manifold $\gamma_u$ have the
obvious counterpart for the stable manifold $\gamma_s$ that can be
readily obtained via reversibility, see Remark \ref{rem:rev}.
\end{rem}

\section{A narrower cone field}\label{sec:narrow}
 
Here our goal is to estimate the angle between stable and unstable manifolds.

More precisely, we wish to prove that there exists two constants $K_+,K_-
\in\R^+$ such that the cone field
$\Co_*(\xi):=\{(1,u)\in\R^2\;|\;K_-(|x|+\sqrt{|y|})\leq u\leq
K_+(|x|+\sqrt{|y|})\}$ contains the unstable direction (by
reversibility we can also define 
the stable cone field $\Co^-_*$). 

\begin{prop}\label{prop:dista}
For each $\xi\in\To^2$ holds true $E^u(\xi)\in\Co_*(\xi)$,
$E^s(\xi)\in\Co_*^-(\xi)$.
\end{prop}

The rest of the section is devoted to the proof of Proposition
\ref{prop:dista}.

Clearly a problem arises only in a neighborhood of zero. Accordingly
the first step is to gain a better understanding of the dynamics near
zero.

\subsection{Near fixed point dynamics}\label{sub:nearfix}
For each $\delta>0$ let $Q_\delta:=[-\delta,\delta]^2$ be a square
neighborhood of zero. The manifolds of the fixed point divide
$Q_\delta$ into four sectors: two thin and two fat (see Figure
\ref{fig:man}). We will discuss explicitly the dynamics in the two
sectors below the unstable manifold (the other two being identical by
symmetry).

\begin{lem}\label{lem:zdyn1}
For each $(x,y)\in\To^2$, $y\leq \gamma_u(x)$, let
$(x_n,y_n):=T^n(x,y)$, then it holds true 
\[
x_n\leq f_u^n(x)\quad \forall n\geq 0.
\]
\end{lem}
\begin{proof}
First note that the trajectory will always remain below the unstable
manifold. Hence, by induction,
\[
\begin{split}
x_{n+1}&=x_n+h(x_n)+y_n\leq x_n+h(x_n)+\gamma_u(x_n)\\
&\leq f_u^n(x)+
h(f_u^n(x))+\gamma_u( f_u^n(x))=f_u^{n+1}(x). 
\end{split}
\]
\end{proof}
The above lemma will suffice to control the dynamics in the thin
sector, more work is needed for the fat one. In fact, when the
trajectories are close to the stable or the unstable manifolds the
above result can still be used (possibly remembering reversibility). On the
other hand when the trajectory is close enough to zero its behavior
is drastically different from the one on the invariant manifolds.

To define more precisely the meaning of ``close to zero'' let us
introduce the parabolic sector $P_M:=\{(x,y)\in Q_1\;|\;|y|\geq Mx^2\}$.
We consider a backward trajectory starting from $x\leq 0$, $y\leq
\gamma_u(x)$, the other possibilities follow by reversibility. Let, as usual, $(x_{n},
y_n):=T^n(x,y)$, $n\in\Z$. Let  $m_+$ be the smallest integer for which
$(x_{-n},y_{-n})\in P_M$, $m$ the largest integer such that
$x_{-m}\leq 0$, and $m_-$ the largest integer for which
$(x_{-n},y_{-n})\in P_M$. Define, see \eqref{eq:hamiltonian}, 
\[
E:=H(x_{-m},y_{-m}).
\]
In addition, define the function $\Upsilon_E: [-1,1]\to\R^-$,
by\footnote{Computing for $y\leq0$ yields
\[
\begin{split}
\Upsilon_E(x)&=-\frac {h(x)}2-\sqrt{\frac{h(x)^2}4+2[G(x)+E-\frac
1{12}h(x)^2](1-\frac 16 h'(x))}\\
&=-\sqrt{2(E+G(x))}(1+\Or(x^2))+\Or(x^3).
\end{split}
\]
}
\begin{equation}\label{eq:upsi}
H(x,\Upsilon_E(x))=E.
\end{equation}
Then, by \eqref{eq:appham}, and since $|y_{-n}|\geq \cs x_{-n}^2$,
\[
\begin{split}
H(x_{-n},y_{-n})-&H(x_{-n},\Upsilon_E(x_{-n}))=H(x_{-n},y_{-n})-H(x_{-m},y_{-m})\leq
\cs \sum_{k=n}^m y_{-k}^4\\
&\leq\cs \sum_{k=n}^m |y_{-k}^3|(x_{-k-1}-x_{-k})
\leq  \cs|y_{-n}^3| |x_{-n}|. 
\end{split}
\]
Accordingly, for $n\leq m$ it holds true
\begin{equation}\label{eq:dyny}
|y_{-n}-\Upsilon_E(x_{-n})|\leq \cs|y_{-n}|^2|x_{-n}|.
\end{equation}
\begin{lem}\label{lem:zdyn2}
In the above described situation, setting $M=\sqrt b$,  the following holds true
\begin{enumerate}
\item $f_u^{-n}(x)\leq x_{-n}\quad \forall\; n\leq m$;
\item $f_s^{k}(x_{-n})\geq x_{-n+k}\quad \forall\;  n\geq m_-,\ k\leq
  n-m_-$;
\item $\sqrt{E}\leq |y_{-n}|\leq 3\sqrt{E}$ for all
  $n\in\{m_+,\dots, m_-\}$;
\item $m_+\leq 2A(Eb^{-1})^{-\frac 14}$;
\item $ 2(12Eb)^{-\frac 14}\leq m_--m_+\leq 4(Eb)^{-\frac 14}$.
\end{enumerate}
\end{lem}
\begin{proof}
The first fact is proven as in Lemma \ref{lem:zdyn1}, the second
follows by reversibility. Hence, by
the results of section \ref{sec:fix}, for $n\leq m$, it follows 
\begin{equation}\label{eq:nest}
|x_{-n}|\leq |f_u^{-n}(x)|\leq \frac{2A|x|}{|x|n+A}.
\end{equation}
Next we want to determine the points $x_{m_+}$ and $x_{m_-}$. The idea
is to use $\eqref{eq:dyny}$ that determines with good precision the
geometry of the trajectories. Let $\bar x$ be defined by
$\Upsilon_E(\bar x)=-M\bar x^2$. Then
\[
\bar x=-\left[\frac{2E}{M^2-\frac b2}\right]^{\frac 14}+\Or(E^{\frac
  54}).
\]  
On the other hand, since by definition $|y_{-m_+}|\geq Mx_{-m_+}^2$
and $|y_{-m_++1}|\leq Mx_{-m_++1}^2$, holds
$|y_{-m_+}-Mx_{-m_+}^2|\leq \cs|x_{-m_+}^3|$.
Hence, by \eqref{eq:dyny} it follows 
\[
|x_{-m_+}-\bar x|\leq \cs\frac{x_{-m_+}^3}{\min\{\bar
x,x_{-m_+}\}}\leq \cs \frac{x_{-m_+}^3}{x_{-m_+}-|x_{-m_+}-\bar x|}.
\]
Solving the above inequality yields
\[
|x_{-m_+}-\bar x|\leq\cs x_{-m_+}^2\leq \cs\bar x^2\leq\cs
\sqrt E.
\]
Analogously $|x_{-m_-}+\bar x|\leq \cs \sqrt E$.
\relax From this (3) and (4) easily follows. Finally,
\[
2|\bar x|\geq |x_{m_-}-x_{m_+}|=|\sum_{n=m_+}^{m_-}y_n|\geq
\cs(m_+-m_-)\sqrt E.
\]
Which implies (5).
\end{proof}

We are now ready to refine our knowledge of the stable and unstable
direction. Let us fix $\varrho\in(0,1/2)$.

\subsection{The cone field--Outside $Q_\varrho$}

The general idea is to take the positive cone field $\Co_+$ (which is
invariant and contains the unstable direction) and to push it forward
in order to obtain a narrower cone 
field. First of all outside $Q_{\sqrt\varrho}$ we have (see \eqref{eq:unsdef})
\[
1\geq F(x,u)\geq F(x,0)\geq 2b \varrho,
\]
where we have chosen $\varrho$ small enough.
Hence the cone field $\Co_0=\{(1,u)\;|\;1\geq u\geq 2b\varrho\}$ is invariant
outside $Q_{\sqrt\varrho}$. It remains to understand what
happens in a neighborhood of the origin of order $\sqrt\varrho$. 

Let us define $\bar u(\xi)$ by the equation $u=F(\xi,u)$. An easy
computation shows
\[
\bar u(\xi)= \frac{-h'(x)+\sqrt{h'(x)^2+4h'(x)}}2= \sqrt{3b}|x|+\Or(x^2).
\]
By reversibility we can restrict ourselves to the case $x\geq 0$, in this
case the only possibility to enter the region $Q_{\sqrt\varrho}$ is
via the fourth quadrant. 

Note that $F(\xi,u)\geq u$ provided $0\leq u\leq \bar u(\xi)$.
This means that if $\xi\not\in Q_{\sqrt\varrho}$ but $T\xi\in
Q_{\sqrt\varrho}$, then the lower bound of the cone $D_\xi T^n\Co_0$
does not decreases until $\sqrt{3b}x_i\leq 2b\varrho$, where
$(x_i,y_i):=T^i\xi$. 

Accordingly, the cone field $\Co_0$ is invariant
also in the fourth quadrant, outside the set $Q_{2\sqrt{\frac b3}\varrho}$.
Now consider the cone field $\Co_1(\xi):=\{(1,u)\;|\;2b\varrho\leq
u\leq L \bar u(\xi)\}$, $L=(3b\varrho)^{-\frac 12}$, for $\xi\in
Q_{\sqrt\varrho}\backslash Q_{2\sqrt{\frac b3}\varrho}$. Clearly, if
$|x|\geq \sqrt \varrho$, then $L\bar u(\xi)\geq 1$. Hence, as the
point enters $Q_{\sqrt\varrho}$, the image of $\Co_0$ is contained in
$\Co_1$, moreover we have already seen that the lower bound is
invariant provided $\xi_i\not\in Q_{2\sqrt{\frac b3}\varrho}$. Let us follow the
upper edge, if $u\leq L\bar u(\xi_i)$, then\footnote{Note that this
computation holds for all $\xi_i\in Q_{\sqrt{\varrho}}\backslash P_M$.}
\[
\begin{split}
F&(\xi_i,u)=F(\xi_i,u)-F(\xi_i,\bar u(\xi_i))+\bar u(\xi_i)\\
&\leq \frac{(L-1)\bar u(\xi_i)}{(1+h'(x_i)+L\bar u(\xi_i))(1+h'(x_i)+\bar
u(\xi_i))}+\bar u(\xi_i)\\
&\leq L\bar u(\xi_i)-(L^2-1)\bar u(\xi_i)^2+\Or(x^3_i)\\
&\leq L\bar u(\xi_{i+1})+\sqrt{3b}|y_{i+1}|-(L^2-1)\bar u(\xi_i)^2+\Or(x^3_i)\\
&\leq L\bar u(\xi_{i+1}),
\end{split}
\]
provided $\xi_{i+1}\not\in P_M$, with $M\leq \sqrt{3b}(L^2-1)$, which is
fine provided $\varrho$ is chosen small enough. The above discussion
can be summarized as follows.
\begin{lem}
\label{lem:conefield-one}
There exists $\varrho>0$: For $\xi\not\in Q_{2\sqrt{\frac b3}\varrho}$
the unstable distribution is contained in $\Co_0$. In addition, in the
set $\{\xi=(x,y)\in Q_{\sqrt{\varrho}}\setminus (P_M\cup
Q_{2\sqrt{\frac b3}\varrho})\;:\; xy\leq 0\}$ the unstable direction
is contained in $\Co_1$.
\end{lem}

To conclude we need to study what happens in a neighborhood of the
origin of order $\varrho$. It is necessary to 
distinguish two possibilities: one can enter below the stable
manifold, and hence be confined in the fat sector, or one can enter
above the stable manifold, thereby being bound to the thin sector.
We will start with the easy case: the second.
 
\subsection{The cone field--Thin sector}

If $x>0$, as soon as the trajectory, at some time $n$, enters
$Q_{2\sqrt{\frac b3}\varrho}$ we have that 
$(1,u)\in\Co_0$ implies $u\geq 2b\varrho\geq \gamma_u'(x_n)$.
Let us consider in $Q_{2\sqrt{\frac b3}\varrho}$ the cone field
$\Co_2:=\{(1,u)\;|\; L\bar u(\xi)\geq u(x)\geq \gamma_u'(x)\}$. Note
that, upon entering in $Q_{2\sqrt{\frac b3}\varrho}$ such a cone contains 
$\Co_1$.\footnote{Note that, in such a case, the trajectory cannot enter in $P_M$.}
Now
\[
F(x,u)\geq F(x,\gamma_u'(x))=\gamma_u'(x+h(x)+\gamma_u(x))\geq
\gamma_u'(x+h(x)+y), 
\]
where we have use that $\gamma_u''(x)\geq 0$ for $x\in[0,\varrho]$,
provided $\varrho$ has been chosen small enough.
\begin{lem}
\label{lem:conefield-two}
In the region $Q_{2\sqrt{\frac b3}\varrho}\setminus(
P_M\cup\{\xi=(x,y)\in\To^2\;:\; y\geq \gamma^u(x) \text{ for }x>0;
y\leq \gamma^u(x) \text{ for }x<0\})$ the
unstable direction is contained in the cone field $\Co_2$. 
\end{lem}
Note that the above lemma suffices for trajectories in the thin sector.
The situation it is not so simple in the fat sector since the lower
bound would deteriorate to zero. A more detailed analysis is needed.

For each $\xi=(x,y)\in\To^2$, for which the unstable direction is defined,
let $(1,u(\xi))$ be the vector in the unstable direction. Define then
$\lambda_{u,n}(\xi, u)$ and $F_n(\xi,u)$ as in formulae \eqref{eq:unsdef}
and \eqref{eq:unstreg} and similarly define the stable
quantities. That is
\begin{equation}\label{eq:constex}
\begin{split}
&D_{T^{-n}\xi} T^{n}(1,u)=:\lambda_{u,n}(\xi,u)(1,F_n(\xi,u))\\
&D_{\xi} T^{-n}(1,-v)=:\mu_{s,n}(\xi,v)(1,-F^-_n(\xi,v))
\end{split}
\end{equation}

\subsection{The cone field--Fat sector}
First of all notice that the trajectory can enter
$Q_{2\sqrt {\frac b3} \varrho}$ either in $P_M$ or outside.
Since the cone field $\Co_2$ for $x\geq
0$, $\xi\not\in P_M$ contains the unstable vector (Lemma
\ref{lem:conefield-two}), we have a good control on the 
unstable vector in both cases until we enter in $P_M$. 
Upon entering $P_M$, we will obtain a very sharp control on the
evolution of the edges of the cone. Let $\xi\not\in P_M$, $T\xi\in
P_M$, and let $\ell_+-1>0$ be the smallest integer such that
$\xi_n\not\in P_M$. By equation \eqref{eq:unsdef}, we have
\[
u_n:=F_n(\xi_n,u)=\sum_{i=1}^{n}
\lambda_{u,i}(\xi_{n-i},u_{n-i})^{-1}h'(x_{n-i})+\lambda_{u,n}(\xi_n,u)^{-1}
u. 
\]
Then, for each $n< \ell_+$, holds true
\[
u_n\leq \sum_{i=1}^n h'(x_{n-i})+u\leq \frac {\cs }{M} 
\left|\sum_{i=1}^n y_i\right|+u. 
\]
By Lemma \ref{lem:zdyn2}-(3),(5), it follows that
we have, for $u\in\Co_2(\xi)$,
\[
u_n\leq \cs_+\sqrt{|y_n|}.
\]
Moreover, remembering \eqref{eq:unsdef} and that  $u\in\Co_2(\xi)$, yields
\[
u_n\geq e^{-2n\cs_+\sqrt{|y_n|}}u\geq \cs u\geq \cs_-\sqrt{|y_n|}. 
\]
Consequently, if for $\xi=(x,y)$ we define the cone
$\Co_3(\xi)=\{\cs_-\sqrt{|y|}\leq u\leq \cs_+\sqrt{|y|}\}$., then the
above results can be written as follows.
\begin{lem}
\label{lem:conefield-three}In $P_M$ the unstable direction
is contained in the cone field $\Co_3$.
\end{lem}
Finally we have to follow the trajectory outside $P_M$ until it exits
from $Q_{2\sqrt{\frac b3}\varrho}$. The upper bound can be treated as
before. Not so for the lower bound.

Let $\xi=(x,y)$ be a point in the
fat sector, $x\leq 0$, $x_{-1}\geq 0$. Then, remembering subsection
\ref{sub:nearfix}, let $E:=H(x,y)$, $u_0=0$ and $u_{n+1}:=F(x_n,
u_n)$. Clearly, $D_{(x,y)}T^n\Co_+\subset\{(1,u)\in\R^2\;|\; u\geq
u_n\}\cup\{0\}$.
\begin{lem}
\label{lem:flowdir}
In the situation described above, for each $n\in\N$, holds true
\[
F(x_n,\Upsilon_E'(x_n))-\Upsilon_E'(x_{n+1})=\Or(|y_n|^{3/2}).
\]
\end{lem}
\begin{proof}
Notice that, since the trajectory lies below the unstable manifold, 
$|y|\geq \cs x^2$. It is then convenient to keep track of the orders
of magnitude only in terms of powers of $y$.
\[
F(x_n,\Upsilon_E'(x_n))=\Upsilon_E'(x_n)-\Upsilon_E'(x_n)^2
+h'(x_n)+\Or(|y_n|^{3/2}).
\]
On the other hand, differentiating \eqref{eq:upsi}, one gets 
\[
\Upsilon_E'(x)=\frac{h(x)}{\Upsilon_E(x)}
-\frac{h(x)^2}{2\Upsilon_E(x)^2} -\frac{h'(x)}2+\Or(|\Upsilon_E|^{3/2}). 
\]
Accordingly, by \eqref{eq:dyny}, 
\[
\begin{split}
\Upsilon_E'(x_{n+1})&=\frac{h(x_n)+h'(x_n)\Upsilon_E(x_n)}{\Upsilon_E(x_n)+h(x_n)}
-\frac{h(x_n)^2}{2\Upsilon_E(x_n)^2}-\frac{h'(x_n)}2 
+\Or(|y_n|^{3/2})\\
&=\Upsilon_E'(x_{n})-\Upsilon_E'(x_{n})^2 +h'(x_n)+\Or(|y_n|^{3/2}),
\end{split}
\]
from which the Lemma easily follows.
\end{proof}
Since $F$ is a contraction in $u$, 
we can estimate
\[
\begin{split}
|u_n-\Upsilon_E'(x_n)|&=|F(x_{n-1}, u_{n-1})-F(x_{n-1},\Upsilon_E'(x_{n-1}))|
+\Or(|y_{n-1}|^{3/2})\\
&\leq |u_{n-1}-\Upsilon_E'(x_{n-1})|+\Or(|y_{n-1}|^{3/2})\\
&\leq|\Upsilon_E'(x_0)|+\Or\left(\sum_{k=0}^{n-1}|y_k|^{3/2}\right)\\
&=\Or\left(|y_0|+\sum_{k=0}^{n-1}\sqrt{|y_k|}(x_{k}-x_{k+1})\right)
=\Or(|y_n|).
\end{split}
\]
We have thus proved that there exists a constant $\cs_0>0$ such that
\begin{equation}\label{eq:fatcone}
u_n\geq \Upsilon_E'(x_n)-\cs_0\Upsilon_E(x_n).
\end{equation}
Hence outside $P_M$ the image of the cone will belong to the
cone field $\Co_3:=\{(1,u)\in\R^2\;|\; u(\xi)\geq
\Upsilon'_{E(\xi)}(x)-\cs_0 \Upsilon_{E(\xi)}(x)\}$. Note that, upon exiting $P_M$,
$\Upsilon_E'(x_n)-\cs_0\Upsilon_E(x_n)\geq \cs_-'\sqrt{|y_n|}$,
provided $\varrho$ is chosen small enough.
The Proposition follows by choosing $\varrho$ small enough and
remembering Lemmata \ref{lem:conefield-one}, \ref{lem:conefield-two}
and \ref{lem:conefield-three}.

\section{An a priori expansion bound}
\label{sec:expansion}

The results of the previous section allow to obtain the following nice
estimate on the expansion in the system.
\begin{lem}\label{lem:expansion}
There exists $K>0$ such that, for each
$\xi=(x,y)\in\To^2\backslash \{0\}$, $n\in\N$ and
$(1,u)\in\Co_*(T^{-n}\xi)=\Co_*(\xi_{-n})$, holds true 
\[
\lambda_{u,n}(\xi,u)\geq e^{-K|x|} \left(K^{-1}
|x|n+1\right)^2 . 
\]
\end{lem}
\begin{proof}
Let us fix $\delta>0$. On the one hand, if the trajectory lies outside
of $Q_\delta$, then 
we have an exponential expansion, on the other hand, if the backward
trajectory enjoys $|x_{-n}|\geq |x_0|$, then equation \eqref{eq:unsdef}
and Proposition \ref{prop:dista} imply
\begin{equation}
\label{eq:trivialexp}
\lambda_{u,n}(\xi,u)\geq (1+K_-|x_0|)^n\geq  e^{-K|x_0|} \left(K^{-1}
|x_0|n+1\right)^2 .
\end{equation}
We say that  the backward orbit of $\xi$ (up to time $n$)  passes $p$ times 
thru
$Q_\delta$ if $\{0\leq k\leq n\;:\;\xi_{-k}\in Q_\delta\}$ consists of
$p$ intervals. The Lemma holds for orbits that pass zero-times thru
$Q_\delta$. Suppose it holds for orbits that pass $p$ times. Let
$\xi_{-n}\in Q_\delta$ and let $m<n$ be the last time $\xi_{-m}\not\in
Q_\delta$ but it passed already $p$ times in $Q_\delta$. Moreover,
suppose that the Lemma holds in $Q_{2\delta}$. Accordingly, for each
$n\geq l$ such that $\xi_{-n}\in Q_\delta$ holds
\[
\begin{split}
\lambda_{u,n}(\xi,u)&\geq e^{-K|x_0|} \left(K^{-1}
|x_0|m+1\right)^2e^{-2K\delta} \left(2K^{-1}
\delta(n-m)+1\right)^2\\
&\geq e^{-K|x_0|} \left(K^{-1}
|x_0|m+1\right)^2\left(K^{-1}\delta(n-m)+1\right)^2\\
&\geq e^{-K|x_0|} \left(K^{-1}|x_0|n+1\right)^2,
\end{split}
\]
provided $\delta$ has been chosen small enough and since it must be $n-m\geq \cs
\delta^{-1}$.  
Thus to prove the Lemma it suffices to prove it for the pieces of
trajectories in $Q_\delta$. There are two cases: a trajectory enters in the thin
sector or in the fat one. Let us consider the thin sector first. 
Set $u_{-j}:=F_{n-j}(\xi_{-j},u)$.
By the usual distortion estimates follows
\[
\begin{split}
\lambda_{u,n}(\xi,u)&=\prod_{j=1}^n(1+h'(x_{-j})+u_{-j})\geq
\prod_{j=1}^n(1+h'(x_{-j})+\gamma_u'(x_{-j}))\\
&\geq
\prod_{j=1}^n(1+h'(f_u^{-j}(x))+\gamma_u'(f_u^{-j}(x)))
=\prod_{j=1}^n(f_u)'(f_u^{-j}(x))\\
&\geq\prod_{j=1}^n\frac
{|f_u^{-j+1}(x)-f_u^{-j}(x)|}{|f_u^{-j}(x)-f_u^{-j-1}(x)|}
e^{-\cs|f_u^{-j+1}(x)-f_u^{-j}(x)|}\\
&\geq e^{-\cs|x_0|}\frac{|x_0|^2}{|f_u^{-n}(x_0)|^2}.
\end{split}
\]
Now, notice that $f_u^{-1}(x)\leq \frac{x}{1+\cs x}$, hence
$f_u^{-n}(x_0)\leq \frac{|x_0|}{1+\cs n|x_0|}$. Thus,
\begin{equation}
\label{eq:bah}
\lambda_{u,n}(\xi,u)\geq e^{-\cs|x_0|}(1+n\cs|x_0|)^2.
\end{equation}
\relax For the fat sector we need only to consider the cases in which $x_0\not\in
Q_\delta$ and $x_0\in Q_\delta$, $x_0\leq 0$ since if $x_0>0$ the
backward trajectory increases the $x$ coordinate. In such cases we
have\footnote{Again, $E$ is chosen to be the energy associated
to the point of the orbit closer to the origin.} 
\[
\begin{split}
\lambda_{u,n}(\xi,u)&=\prod_{j=1}^n(1+h'(x_{-j})+u_{-j})\geq
\prod_{j=1}^n(1+\Upsilon_E'(x_{-j})-\cs_0\Upsilon_E(x_{-j}))\\
&\geq e^{\sum_{j=1}^n\Upsilon_E'(x_{-j})-2\cs_0\Upsilon_E(x_{-j})}\\
&\geq e^{-\sum_{j=1}^n\frac{\Upsilon_E'(x_{-j})}{\Upsilon_E(x_{-j})}
(x_{-j-1}-x_{-j}) -3\cs_0\sum_{j=1}^n(x_{-j-1}-x_{-j})} \\
&\geq e^{-\cs|x_0|}e^{-\int_{x_0}^{x_{-n}}\frac{\Upsilon_E'(z)}{\Upsilon_E(z)}dz}
= e^{-\cs|x_0|}\frac{\Upsilon_E(x_0)}{\Upsilon_E(x_{-n})}.
\end{split}
\]

Let $n_*\in\N$ be the last integer for which $|G(x_{-n})|\geq
E$, then for $n\leq n_*$ we have 
\[
\lambda_{u,n}(\xi)\geq e^{-\cs|x_0|}\sqrt{\frac
{E+G(x_0)}{E+G(x_{-n})}} 
\geq e^{-\cs|x_0|}\sqrt{\frac
{G(x_{-n})+G(x_0)}{2G(x_{-n})}}.
\]
On the other hand comparing the backward motion with the
backward motion on the stable manifold, as we did before with the
unstable,\footnote{Here we use the inequality  
\[
\sqrt{\frac{1+(1+a)^4}2}\geq (1+\frac a2)^2.
\]
}
\[
\lambda_{u,n}(\xi,u)\geq
e^{-\cs|x_0|}\sqrt{\frac{1+(n|x_0|\cs+1)^4}2}
\geq e^{-\cs|x_0|}(1+n|x_0|\cs)^2.
\]
Next, let us consider $n\in\{n_*,\dots, m\}$, where $m$ is the larger
integer such that $x_{-m}\leq 0$, we have $2\sqrt E\geq 
\Upsilon_E(x_{-n})\geq \sqrt {2E}$. 
\[
\begin{split}
\lambda_{u,n-n_*}(\xi,u)
&\geq e^{-\cs_2|x_{n_*}|}
e^{-\int_{x_{-n_*}}^{x_{-n}}\frac{\Upsilon_E'(z)}{\Upsilon_E(z)}dz} 
\geq e^{-\cs|x_0|}\frac{\Upsilon_E(x_{-n_*})}
{\Upsilon_E(x_{-n})}\\
&\geq  e^{-\cs|x_0|}\sqrt{\frac 32}\geq e^{-\cs|x_0|}(1+\cs
|x_{-n_*}|(n-n_*))^2, 
\end{split}
\]
where, in the last line, we used Lemma \ref{lem:zdyn2}-(5).
By symmetry it will be enough to wait another time $m$ to have
$|x_{-2m}|\geq \frac 12|x_0|$, after which the expansion is assured by
the estimate \eqref{eq:trivialexp}.
\end{proof}

Next we need to have similar estimates for the stable contraction. By
\eqref{eq:constex} 
\begin{equation}\label{eq:stabdef}
\begin{split}
\mu_s(v)&:=1+v=\mu_{s,1}(\xi,v)\\
F^-_1(\xi,v)&=h'(x-y)+\frac v{1+v}.
\end{split}
\end{equation}
It is immediate to check that
$D_{(x,y)}T^{-1}(1,-v)=\mu_s(v)(1,-F^-_1((x,y),v))$ and
$\mu_{s,n}(\xi,v_{0}):=\prod_{i=0}^{n-1}\mu_s(v_{-i})$, where 
$v_0=v$ and $v_{-i-1}:=F^-_1(T^{-i}\xi,v_{-i})$.

An interesting way to transform information on expansion into information on
contraction is to use area preserving. 
\begin{lem}\label{lem:areap}
Let $\xi\in\To^2$, then for each $n\in\N$, $u,v\geq 0$ let $u_{-n}=u$,
$v_0=v$, $\xi_{-n}=T^{-n}\xi$ and $u_{-k+1}=F(\xi_{-k+1}, u_{-k})$,
$v_{-k-1}=F^-(\xi_{-k},v_{-k})$. Then
\[
\mu_{s,n}(\xi, v_0)(v_{-n}+u_{-n})=\lambda_{u,n}(\xi, u_{-n})(u_0+v_0).
\]
\end{lem}
\begin{proof}
Calling $\omega$ the standard symplectic form we have
\[
\begin{split}
\mu_{s,n}(\xi,v_0)(v_{-n}+u_{-n})&=\omega(D_{\xi}T^{-n}(1,-v_0),(1,u_{-n}))\\
&= \omega((1,-v_0), D_{\xi_{-n}} T^n(1,u_{-n}))=\lambda_{u,n}(\xi,
u_{-n})(u_0+v_0). 
\end{split}
\]
\end{proof}
%\begin{changebar}
The following is an immediate corollary of Lemmata \ref{lem:areap} and
\ref{lem:expansion}. 
\begin{cor}\label{lem:contraction}
For each $\xi=(x,y)\in \To^2$ and $n\in\N$  holds
\[
\mu_{s,n}(\xi, v_0)\geq e^{-\cs|x_{0}|} (\cs^{-1}|x_{0}|n+1)^2\frac
{u_0+v_0}{u_{-n}+v_{-n}} \quad\forall n\in\N.
\]
\end{cor}
%\end{changebar}

All the other expansion estimates can be obtained by reversibility.

\section{Distributions--regularity}
\label{sec:regularity}

Let $(1,u(\xi)),\,(1,-v(\xi))$ be the unstable and stable directions,
respectively. We will then use the short hand
$\lambda_{u,n}(\xi):=\lambda_{u,n}(\xi,u(\xi_{-n}))$ and
$\mu_{s,n}(\xi):=\mu_{s,n}(\xi,v(\xi_{n}))$. 
\begin{lem}\label{lem:distreg0}
The unstable distribution is continuous in $\To^2$.
\end{lem}  
\begin{proof}
Notice that, for $\xi=(x,y)$, $\xi_n:=T^n\xi$, iterating formula
\eqref{eq:unsdef}, in analogy with \eqref{eq:unstreg}, holds true
\begin{equation}\label{eq:unstreg-bis}
\begin{split}
u(x)-u(z)=&\lambda_{u,n}(x)^{-1}\lambda_{u,n}(z)^{-1}(u(x_{-n})-u(z_{-n}))\\
&\ \ +\sum_{k=1}^{n}\lambda_{u,k}(x)^{-1}\lambda_{u,k}(z)^{-1}
(h'(x_{-k})-h'(z_{-k}))
\end{split}
\end{equation}
By Lemma \ref{lem:expansion}, we can take the limit $n\to\infty$ in
the above formula provided  $x\neq 0$, and obtain a
uniformly convergent series from which the continuity follows.
If $\xi\neq 0$ then $x_{-1}\neq 0$ and \eqref{eq:unsdef} implies
\begin{equation}
\label{eq:uma}
u(\xi)=\lambda_{u}(\xi_{-1})^{-1}h'(x_{-1})
+\lambda_{u}(\xi_{-1})^{-1}u(\xi_{-1}) ,
\end{equation}
hence the continuity at $\xi\neq 0$ follows. We are left with the
continuity at the origin, but this is already implied by Proposition
\ref{prop:dista}. 
\end{proof}

This means that we can extend the invariant unstable distribution
(that, up to now, where defined--by Pesin theory--only almost everywhere) to a
continuous everywhere defined vector field. The same statement holds
for the stable vectors by reversibility.

Given a continuous vector field there exists integral curves. Since
we do not know yet if the vector fields are Lipschitz, it does not
follows automatically that from a given point there exits only one
integral curve, yet this follows by standard dynamical
arguments. Clearly such integral curves are nothing else than the
stable and unstable manifolds that are therefore everywhere
defined. In addition, remember that, by general hyperbolic theory, the
foliations are absolutely continuous, it follows that the above
everywhere defined foliations are continuous. Unfortunately, for the
following much sharper regularity information is  needed, this is
obtained in the rest of the section.

Let us call $\partial^u, \partial^s$ the derivative along the unstable
and the stable vector fields, respectively. 

\begin{lem}\label{lem:uderbound}
The vector field $u$ is $\Co^1$ along the unstable manifolds, apart from the
origin, moreover
\[
|\partial^u
u(\xi)|\leq C\quad
\forall \xi\neq 0. 
\]
\end{lem}
\begin{proof}
If $\xi$ is outside of a 
neighborhood of the origin of size $\delta$, then by Lemma
\ref{lem:expansion}, \eqref{eq:unstreg-bis} we have, in analogy with
the arguments leading to \eqref{eq:unstder},
\begin{equation}
\label{eq:uder-u-bond}
|\partial^u
u(\xi)|=\left|\sum_{k=1}^{\infty}\lambda_{u,k}(x)^{-3}h''(x_k)\right| \leq \cs
\sum_{n=0}^{\infty}(\delta n+1)^{-6}\leq \cs. 
\end{equation}
Since the series converges uniformly the $\Co^1$ property follows. To
obtain a uniform bound more work is needed.
If $|\xi|<\delta$, formula \eqref{eq:uma} implies
\[
|\partial^u u(\xi)|\leq\lambda_u(\xi_{-1})^{-3}|h''(x_{-1})|+
\lambda_u(\xi_{-1})^{-3} |\partial^u u(\xi_{-1})|=:\Psi(\xi_{-1},
|\partial^u u(\xi_{-1})|).
\]
A simple computation, remembering Proposition \ref{prop:dista}, shows
that 
\[
\Psi(\xi,\varrho)\leq7b|x|+\frac{\varrho}{1+3|u(\xi)|}\leq
7b|x|+\frac{\varrho}{1+3K_-|x|} \leq \varrho ,
\]
provided $\varrho\geq \frac{7b(1+3K_-)}{3K_-}$. Accordingly, for $\rho$
large enough, we have 
$|\partial^u u(\xi)|\leq \rho$, for all $\xi$. 
\end{proof}

It remains to investigate the regularity of the unstable distribution
along the stable direction.

\begin{lem}\label{lem:distreg1}
The unstable distributions are $\Co^1$ along stable manifolds, apart
from the origin. Moreover
\[
|\partial^s u(\xi)|\leq \cs \quad \forall \xi\neq 0.
\]
\end{lem}  
\begin{proof}
Let us fix some arbitrary neighborhood of the origin.
Let $x,z\in W^s$ outside such a neighborhood. Let $W^s_0$ be the
piece of stable manifold between such two points. Clearly
$W^s_{n}:=T^nW^s_0$ grows for negative $n$. Let $n(x,z)$ be the largest
integer for which $|W^s_{-n}|\leq |W^s_0|^{\frac 14}$. Our first
result is a distortion bound.
\begin{sublem}\label{slem:dist}
For each $n\leq n(x,z)$ and $\xi\in
W^s_{0}$, holds
\[
\cs^{-1}\frac{|W^s_{-n}|}{|W^s_{0}|}\leq \mu_{s,n}(\xi)\leq
\cs\frac{|W^s_{-n}|}{|W^s_{0}|} .
\]
\end{sublem}
\begin{proof}
If the backward orbit spends at least half of the time outside the
neighborhood, then $W^s_{-n}$
grows exponentially fast, hence $n(x,z)\leq \cs\ln|W^s_0|^{-1}$ and
$\sum_{i=0}^{n(x,z)}|W^s_{-i}|\leq \cs$. If this is not the case,
the worst possible situation is when $W^s_{-m}$ is the closest
to the origin and all the trajectory lies in the neighborhood. In such a
case, letting $m:=n(x,z)$, 
\[
|W^s_{m}|=\int_{W^s_0}\mu_{s,m}(z)dz\geq
\int_{W^s_0}\frac{\theta(z)(\cs^{-1}|z|m+1)^2}{2\theta(T^{-m}z)} dz,
\]
where $\theta(\zeta)=u(\zeta)+v(\zeta)$ is the separation between the
stable and the unstable directions at the point $\zeta$ and we have
used Lemma \ref{lem:contraction}. Now  Proposition \ref{prop:dista} and
Lemma \ref{lem:zdyn2}-(1) imply $\theta(T^{-m}z)\leq \cs m^{-1}$
outside the parabolic sector, while Lemma \ref{lem:zdyn2}-(3,4,5) show
that the same estimates remain in $P_M$ as well. Accordingly,
\[
|W^s_{-m}|\geq \cs m^3|W^s_0|.
\]
That is $m\leq \cs |W^s_0|^{-\frac 14}$, and
\[
\sum_{i=0}^{n(x,z)}|W^s_{-i}|\leq n(x,z)|W^s_0|^{\frac 14}\leq \cs.
\]
The above estimate readily implies that, for each $\xi,\eta\in W^0_s$,
\[
e^{-\cs |W^s_{-i}|}
\leq\frac{\mu_s(\xi_{-i},v(\xi_{-i}))}{\mu_s(\eta_{-i},v(\eta_{-i}))}
\leq e^{\cs |W^s_{-i}|},
\]
where we have used Lemma \ref{lem:uderbound} for the stable manifold.
Accordingly,
\[
e^{-\cs\sum_{i=0}^m |W^s_{-i}|}
\leq\frac{\mu_{s,n}(\xi)}{\mu_{s,n}(\eta)}
\leq e^{\cs \sum_{i=0}^m|W^s_{-i}|},
\]
from which the Lemma readily follows.
\end{proof}
By Lemma \ref{lem:areap} it follows, letting again $m:=n(x,z)$,
\[
\begin{split}
\lambda_{u,m}(x)^{-1}\lambda_{u,m}(z)^{-1}
&=\lambda_{u,m}(x)^{-1}\sqrt{\lambda_{u,m}(z)^{-2}}\\
&\leq
\cs\left(\theta(x_{-m})\mu_{s,m}(x)
\sqrt{\lambda_{u,m}(z)\theta(z_{-m})\mu_{s,m}(z)}\right)^{-1}
\end{split}
\]
As before the worst case is clearly when $W^s_{-m}$ is the closest to
the origin.
In such a case, consider that at least one of the two end points of
$W^s_{-n(x,z)}$ must be at a distance $\cs |W^s_{-n(x,z)}|$ from the
fixed point, let us say $T^{-n(x,z)}z$, hence
$\theta(T^{-n(x,z)}z)\geq \cs|W^s_{-n(x,z)}|$. In addition,
$\theta(x_{-m})\geq \cs m^{-1}$. Indeed, this follows from Lemma
\ref{lem:zdyn2}-(3,4,5) if the trajectory ends in $P_M$. If the
trajectory lies outside $P_M$ then it approaches the origin slower than
the dynamics $x-\cs x^2$, which implies $x_{-m}\geq \cs
m^{-1}$. Furthermore by using the above facts, Lemma \ref{lem:expansion}, Sub-lemma
\ref{slem:dist} and the definition of the stopping time $m$ yields
\[
\lambda_{u,m}(x)^{-1}\lambda_{u,m}(z)^{-1}\leq \cs|x-z|^{\frac
34}m\frac 1{m} \sqrt{\frac{|x-z|^{\frac 
34}}{|x-z|^{\frac 14}}} \leq \cs|x-z|.
\]
Since we know that $u$ is a uniformly continuous function it follows
\[
\lim_{z\to
x}\lambda_{u,m}(x)^{-1}\lambda_{u,m}(z)^{-1}
\frac{|u(T^{-m}x)-u(T^{-m}z)|} 
{|x-z|}=0.
\]
Accordingly, by formula  \eqref{eq:unstreg}, 
\[
\begin{split}
u'(x)=&\lim_{z\to x}
\sum_{n=0}^{m}\lambda_{u,n}(x)^{-1}\lambda_{u,n}(z)^{-1}
\frac{h'(T^{-n}x)-h'(T^{-n}z)}{x-z}\\
=& \lim_{z\to x}
\sum_{n=0}^{m}\lambda_{u,n}(x)^{-1}\lambda_{u,n}(z)^{-1}
h''(T^{-n}\zeta_{n}) \mu_{s,n+1}(\zeta_n),
\end{split}
\]
for some $\zeta_n\in W_0$.
But $|h''(T^{-n}\zeta_n)|\leq \cs\theta(T^{-n}z)$, hence
\[
\begin{split}
&\lambda_{u,n}(x)^{-1}\lambda_{u,n}(z)^{-1}\theta(T^{-n}x)
\mu_{s,n+1}(x)\leq \cs\lambda_{u,n}(z)^{-1}\\
&\lambda_{u,n}(x)^{-1}\lambda_{u,n}(z)^{-1}\theta(T^{-n}z)
\mu_{s,n+1}(z)\leq \cs\lambda_{u,n}(x)^{-1}.
\end{split}
\]
Remembering Sub-Lemma \ref{slem:dist} the uniform convergence of the
series follows and yields the formula
\begin{equation}
\label{eq:uder-s-bound}
\partial^s u(x)=\sum_{n=0}^{\infty}\lambda_{u,n}(x)^{-2}\mu_{s,n+1}(x)
h''(T^{-n}x).
\end{equation}
Given the arbitrariness of the neighborhood of zero, the above formula
holds for each $x\neq 0$ and, since the series converges uniformly,
the $\Co^1$ property follows. We can now conclude the Lemma.
By Lemma \ref{lem:areap} follows
\[
\begin{split}
|\partial^s u(x)|&\leq\cs \sum_{n=0}^\infty
\lambda_{u,n}(x)^{-2}\mu_{s,n}(x)\theta(T^{-n}x) \\
&\leq \cs \sum_{n=0}^\infty\lambda_{u,n}(x)^{-1}\theta (x)
\leq \cs \sum_{n=0}^\infty\frac{|x|}{(|x|n+1)^2}\leq \cs .
\end{split} 
\]
\end{proof}
\begin{rem}
Notice that the symmetrical statements follow by reversibility.
\end{rem}
The final result on the regularity of the foliations can be  stated as
follows.
\begin{lem}
\label{lem:c1-reg}
The stable and unstable vector fields are $\Co^1(\To^2\setminus\{0\})$ and, more
precisely, for each $\xi\in\To^2\setminus\{0\}$,
\[
|Du(\xi)|\leq \cs \theta(\xi)^{-1}.
\]
\end{lem}
\begin{proof}
The $\Co^1$ property follows from Lemma 19.1.10 of \cite{KH}.
Then the size of the derivative can be easily estimated by the size of
the partial derivatives in the stable and unstable directions divided by
the angle between them.
\end{proof}
\begin{rem}
In fact, it is likely that with a little more work one can show that
the foliations are $\Co^{\frac 32-\ve}$, but we do not investigate this
possibility since it is not needed in the following.
\end{rem}
\section{Holonomy}
\label{sec:holo}
There exists $\cs_1>0$ such that, given two close by stable manifolds
$W^s_1$, $W^s_2$ we can define the {\em unstable holonomy}
$\Psi^u:W^s_1\to W^s_2$ by $\{\Psi^u(\xi)\}:=W^u(\xi)\cap W^s_2$. 
Let $D_r:=\{\zeta=(z_1,z_2)\in\R^2\;:\;|z_1|\leq r;|z_2|\leq r^2\}$.
\begin{lem}\label{lem:holo}
For each $W^s_1,W^s_2$ disjoint from $D_r$ and $\xi\in W^s_1$, holds
\[
|J\Psi^u(\xi)-1|\leq \cs r^{-1}\|\Psi^u(\xi)-\xi\|.
\]
Provided $\|\Psi^u(\xi)-\xi\|\leq \cs_1 r$.
\end{lem}
\begin{proof}
Let $\gamma_s,\tilde\gamma_s:[-\delta,\delta]\to\R^2$ be $W^s_1,
W^s_2$, respectively, parametrized by arc-length. Also, let
$\Gamma:[-\delta,\delta]^2\to\R^2$, be such that $\Gamma(0,0)=\xi$,
$\Gamma(s,0)=\gamma_s(s)$ 
and $\Gamma(s,t)$ be the unstable manifold, parametrized by
arc-length, of $\Gamma(s,0)$ and, finally,
$\Gamma(0,\rho):=\Psi^u(\xi)$. Note that $\Gamma(s,t)$ can be obtained
integrating the unstable vector field starting from $\Gamma(s,0)$,
hence Lemma \ref{lem:c1-reg} and the standard results on the continuity
with respect to the initial data imply $\Gamma\in\Co^1$. By the
transversality of the stable and unstable manifolds there exist
$\tau,\sigma:[-\delta,\delta]\to\R$ such that
$\Gamma(s,\tau(s))=\Psi^u(\gamma_s(s))=\tilde\gamma(\sigma(s))\in
W^s_2$. Calling $\eta(s)$ the unit vector perpendicular to $\tilde
\gamma'(s)$, by the implicit function theorem, it follows 
\begin{equation}
\label{eq:implicit}
\begin{split}
\tau'(s)&=-\frac{\langle \eta(s),\partial^s\Gamma(s,\tau(s))\rangle}
{\langle\eta(s), \partial_t\Gamma(s,\tau(s))\rangle}\\
\sigma'(s)&=\langle\tilde\gamma_s'(s),\partial^s\Gamma(s,\tau(s))\rangle
-\frac{\langle\tilde 
\gamma'_s(s),\partial_t\Gamma(s,\tau(s))\rangle\;\langle
\eta(s),\partial^s\Gamma(s,\tau(s))\rangle} 
{\langle\eta(s), \partial_t\Gamma(s,\tau(s))\rangle}
\end{split}
\end{equation}
where, clearly, $\sigma'(s)=J\Psi^u(\gamma_s(s))$. Calling
$v^u(\eta)$, $\eta\in\To^2$,
the unit vector in the unstable direction at $\eta$ and $v^s(\eta)$ the
stable one, one has
$\partial_t\Gamma(s,\tau(s))=v^u(\Gamma(s,\tau(s)))$. On the other
hand, setting $V(s,t):=\partial^s\Gamma(s,t)-v^s(\Gamma(s,t))$, holds
$V(s,t)=0$ for $t=0$, but for $t\neq 0$, in general, it will be
$V(s,t)\neq 0$. Yet, it is possible to estimate it by differentiating
$
\Gamma(s,t)=\Gamma(s,0)+\int_0^tv^u(\Gamma(s,t'))dt'
$
which yields
\[
\partial^s\Gamma(s,t)=v^s(\Gamma(s,0))+\int_0^t Dv^u(\Gamma(s,t'))
\partial^s\Gamma(s,t')dt'.
\]
Lemmata \ref{lem:c1-reg} and \ref{lem:distreg1}
imply that $\|Dv^u\|\leq \cs r^{-1}$ and  $\|Dv^u v^s\|=|\partial^s
v^u|\leq \cs$, hence 
\[
\|V(s,t)\|\leq \cs r^{-1}\int_0^t\|V(s,t')\|dt'+\cs t.
\]
By Gronwal, it follows,  provided $t\leq \cs \rho$ and $\rho\leq \cs_1
r$, for $\cs_1$ small enough,
\begin{equation}
\label{eq:Vest}
\|V(s,t)\|\leq \cs t.
\end{equation}
Accordingly, by the second of \eqref{eq:implicit} and \eqref{eq:Vest},
it follows
\[
|\sigma'(0)-1|\leq \cs r^{-1}\rho.
\]
\vskip-.5cm
\end{proof}

\section{Random perturbations}
\label{sec:random}

The density of a measure with respect to Lebesgue evolves as
\[
\Lp f:=f\circ T^{-1}.
\]
We will then construct a random perturbation by introducing the
convolution operator
\begin{equation}\label{eq:conv}
Q_\ve f(x):=\int_{\To^2}q_\ve(x-y)f(y)dy.
\end{equation}
Where we assume
\begin{enumerate}
\item $q_\ve(\xi):=\ve^{-2}\bar q(\ve^{-1}\xi)$; $\bar q\in\Co^{\infty}(\R^2,\R^+)$;
\item $\int_{\R^2}\bar q(\xi)d\xi=1$;
\item $\bar q(\xi)=0$ for each $\|\xi\|\geq 1$;
\item $\bar q(\xi)=1$ for each $\|\xi\|\leq \frac 12$.
\end{enumerate}
We define then
\begin{equation}\label{eq:random}
\Lp_\ve:=Q_\ve\Lp^{n_\ve},
\end{equation}
where $n_\ve$ will be chosen later.

Notice that 
\begin{equation}\label{eq:ker1}
\Lp_\ve^2f(x)=\int_{\To^4} q_\ve(x-y)q_\ve(T^{-n_\ve}y-T^{n_\ve}z)f(z)
dz dy :=\int_{\To^2}\Ka_\ve(x,z)f(z)dz.
\end{equation}

We have thus a kernel operator that can be investigated with rather
coarse techniques. It turns out to be convenient to define the
associated kernel
\begin{equation}\label{eq:kernel2}
\bar\Ka_\ve(x,z):=\Ka_\ve(x,T^{-n_\ve}z)=\int_{\To^2}q_{\ve}(x-y)
q_\ve(T^{-n_\ve}y-z)m(dy). 
\end{equation}
For further use let us define
\begin{equation}\label{eq:levelset}
\begin{split}
&D_r:=\{z=(z_1,z_2)\in\To^2\;|\;|z_1|\leq r;\; |z_2|\leq r^2\}\\
&B_r(\xi):=\{\eta\in\To^2\;;\;\|\xi-\eta\|<r\}.
\end{split}
\end{equation}
The following is a relevant fact used extensively in the sequel.
\begin{lem}
\label{lem:levelset-escape}
There exists $\cs_3,R>0$ such that, for each $\delta<R^2$, if
$B_\delta(\xi)\subset D_{R/2}$, then there exists $\eta\in
B_{\frac 34\delta}(\xi)$ such that $T^n B_{\frac \delta 4} (\eta)\cap
D_{R}=\emptyset$ for some $n\leq \cs_3 \delta^{-\frac 12}$.
\end{lem}
\begin{proof}
%\begin{changebar}
If $\xi$ belongs to the first or third
quadrant, then $T^n\xi$ is escaping from the origin. 
In such a case, if $B_{\frac 58\delta}(\xi)$ belongs to the thin sector we choose
$\eta\in B_{\frac{\delta}2} (\xi)$, $|y|\geq \frac 38 \delta $.
Clearly, if $(x,y)\in B_{\frac \delta 4}(\eta)$, $\cs x^2\geq |y|\geq
\frac \delta 4$. On the other hand 
$|y_n|\geq n|x|^3$ while $|x_n|\leq |x|+n\cs|x_n|^2$. So, if $n\leq
\cs |x|^{-1}\leq \cs\delta^{-\frac 12}$, we have $|y_n|\geq
\cs|x_n|^2$. After that we can compare the dynamics with one of the
type $x\mapsto x+\cs x^2$, hence after a time at most $\cs
\delta^{-\frac 12}$ the $T^n B_{\frac \delta 4}(\eta)$ will exit $D_R$.
If the above does not apply, then
one can take a ball of radius $\delta/4$ belonging completely to the
fat sector and centered at a point in $B_\delta(\xi)$. Then the results
of subsection \ref{sub:nearfix} easily implies the lemma. 
%\end{changebar}
If, on the contrary, $\xi$ belong to the second or fourth quadrant, then its
trajectory may approach the origin in an arbitrary manner (even
asymptotically, if the point belongs to the stable manifold). In such
a case we can take a point $\eta\in B_{\delta}(\xi)$ at, at least, a vertical
distance $\frac 34\delta$ from the stable manifold and such that
$B_{\frac 34\delta}(\eta) \subset B_\delta(\xi)$. Again from the
results of subsection \ref{sub:nearfix} it follows that such a ball
will exit $D_R$ in a time at most $\cs_3\delta^{-\frac12}$.
\end{proof}

\begin{lem}\label{lem:estim}
There exit constants $\sigma,\cs_3>0$ such that if $n_\ve=\cs_3\ve^{-\frac 12}$ holds
\[
\bar\Ka_\ve(x,z)\geq \sigma\quad \forall x,z\in\To^2.
\]
\end{lem}
\begin{proof}
It is trivial to see that
\[
\bar\Ka_\ve(x,z)\geq \ve^{-4}m(B_{\ve/2}(x)\cap T^{n_\ve}B_{\ve/2}(z)).
\]
Accordingly, by Lemma \ref{lem:levelset-escape} there exists two
balls, of radius $\frac \ve 8$, $\bar B_1\subset B_{\ve/2}(z)$ and
$\bar B_2\subset B_{\ve/2}(x)$ that are outside of $D_{r}$, $r\geq
\sqrt{\frac \ve 2}$, and whose images will be outside
of a neighborhood of the origin or order one in a time less than $\cs
\ve^{-\frac 12}$, forward and backward in time, respectively. Given two 
unstable manifolds in $B_1$ at a distance larger than $\bar c r\ve$, for
some appropriate $\bar c$, then no stable 
manifold will intersect both manifolds inside the ball  $B_1$. We
can thus consider $\cs r^{-1}$ unstable manifolds such that no stable
manifolds intersect two of them in $B_1$. Around each such
manifold we can construct a strip by moving along the stable manifold
by $\cs \ve$. We obtain in this way $\cs r^{-1}$ disjoint strips
each of area $\cs r\ve^2$, whose union covers a fixed fraction of the area of
$B_1$. After a time less that $\cs \ve^{-\frac 12}$ such strips will be outside a
neighborhood of zero, their length may have increase considerable, if so
we will subdivide them into strips of length $\ve$. Since now the
stable and unstable manifold are at a fixed angle and by the usual
distortion arguments, such strips are essentially rectangular. At this point,
by Lemma \ref{lem:expansion}, it will suffice to wait a time
$\ve^{-\frac 12}$ to insure that each such strip will acquire length
at least $\frac 12$ in the unstable direction. We thus iterate for such
a time and, if one strip becomes longer than one, we subdivide it into
pieces of length between $\frac 12$ and one. Finally, fix some box
$\Lambda$ of some fixed size $C$ away from  the origin with sides
approximately parallel either to the stable or to the unstable
directions. 
By mixing it suffices to wait a fixed time to be sure that a fixed
percentage of each one of the above mentioned strips will intersect
the box. In addition, it is possible to insure that such strips cut
the box from one stable side to the other.

We can then write $m(B_{\ve/2}(x)\cap
T^{n_\ve}B_{\ve/2}(z))=m(T^{-n_\ve/2}B_{\ve/2}(x)\cap
T^{n_\ve/2}B_{\ve/2}(z))$ since the same considerations done above for
the unstable manifold can be done, iterating 
backward, for the stable manifold it follows that a fixed percentage of
$T^{-n_\ve/2}B_{\ve/2}(x)$ and a fixed percentage of
$T^{n_\ve/2}B_{\ve/2}(z))$ will intersect $\Lambda$ and hence each
other. In fact each one of the above constructed strips in the unstable
direction will intersects each one of the strips in the stable
direction. By the usual distortion estimates, this implies that the
intersection among any two such strip has a measure proportional to the
product of the measure of the two strips, hence
\[
m(B_{\ve/2}(x)\cap T^{n_\ve}B_{\ve/2}(z))\geq \cs
m(B_{\ve/2}(x))m(B_{\ve/2}(z)),
\]
and the lemma.
\end{proof}

\begin{lem}\label{lem:l1}
For each $f\in L^1$, $\int f=0$ holds
\[
\|\Lp_\ve^n f\|_1\leq (1-\sigma)^{n/2}\|f\|_1.
\]
\end{lem}
\begin{proof}
Note that $\Lp_\ve 1=
\Lp_\ve^*1=1$ and let $\M^+_\ve=\{x\in\To^2\;|\;\Lp_\ve^2 f\geq 0\}$;
$\M_+=\{x\in\To^2\;|\; f\geq 0\}$, 
then, since $\int \Lp_\ve f=\int f=0$,
\[
\begin{split}
\|\Lp_\ve^2 f\|_1&=2\int_{\M^+_\ve}dx\,\Lp_\ve^2 f
=2\int_{\M^+_\ve}dx\int_{\To^2}dy\, \Ka_\ve(x,\,y) 
f(y)\\
&=2\int_{\To^2}dy\, f(y)\int_{\M^+_\ve}dx[ \Ka_\ve(x,\,y)-\sigma]\\
&\leq 2\int_{\M_+}dy\, f(y)\int_{\To^2}dx [\Ka_\ve(x,\,y)-\sigma]
=2(1-\sigma)\int_{\M_+}dy\, f(y)\\
&=(1-\sigma)\|f\|_1.
\end{split}
\]
\end{proof}
Let $v^{u,s}=(v^{u,s}_1,v^{u,s}_2)$ be the unit tangent vector fields in the
unstable and stable direction, respectively. Clearly
$|\partial^u(\Lp^i f)|\leq |\partial^u f|$, while 
$|\partial^s(g\circ T^i)\leq |\partial^s g|$.
\begin{lem}\label{lem:appr}
For each $f,g\in\Co^{1}(\To^2,\R)$ holds
\[
\left|\int Q_\ve f g-\int f g\right|\leq \cs
\ve\{\|f\|_\infty+\|\partial^u f\|_{L^1(\nu)}\}
\{\|g\|_\infty+\|\partial^s g\|_{L^1(\nu)}\}, 
\] 
where $\nu$ is the measure defined by $\nu(h):=\int
d\rho \int_{\partial D_{\rho}}h$.
\end{lem}
\begin{proof}
It is convenient to introduce the following change of
variables. For each $x,y\in\To^2$ close enough, let us call
$[x,y]=W^u_\delta(x)\cap W^s_\delta(y)$, note that by Lemma
\ref{prop:dista} such a point is
always well defined provided $d(x,y)\leq \cs \delta^2$. We consider then
the change of variable $\Phi:\To^2\times \To^2\to\R^2\times\To^2$,
\begin{equation}\label{eq:changv}
\begin{split}
\xi:&=x-y\\
\eta:&=[x,y].
\end{split}
\end{equation}
Due to the absolute continuity of the holonomies the above change of
variable is absolutely continuous. Clearly, $\partial^u_x \eta=0$ since
moving along $W^u(x)$ 
does not change the intersection point with $W^s(y)$. On the other
hand $\partial^s_x \eta=J\Psi^u_x v^s(\eta)$, since moving $x$ along
$W^s(x)$ moves $\eta$ on $W^s(y)$ by an amount determined exactly by the unstable
holonomy $\Psi^u$ between $W^s(x)$ and $W^s(y)$. By similar arguments
and a straightforward computations 
\[
\begin{split}
\partial_{x_1}\eta&=v_2^u(x)\det\begin{pmatrix} v^s(x)
&v^u(x)\end{pmatrix}^{-1}J\Psi^u_x v^s(\eta)\\
\partial_{x_2}\eta&=v_1^u(x)\det\begin{pmatrix} v^u(x)
&v^s(x)\end{pmatrix}^{-1}J\Psi^u_x v^s(\eta)\\
\partial_{y_1}\eta&=v_2^s(y)\det\begin{pmatrix} v^u(y)
&v^s(y)\end{pmatrix}^{-1}J\Psi^s_y v^u(\eta)\\
\partial_{y_2}\eta&=v_1^s(y)\det\begin{pmatrix} v^s(y)
&v^u(y)\end{pmatrix}^{-1}J\Psi^s_y v^u(\eta).
\end{split}
\]
In fact, calling $\theta(x)$ the sine of the angle between stable and
unstable directions at the point $x$ and $v_\perp$ the orthogonal
unit vector to $v$, holds 
\[
\begin{split}
\partial_x\eta&=J\Psi^u_x\theta(x)^{-1}\,  |v^s(y)\rangle\langle v^u_\perp(x)|\\
\partial_y\eta&=J\Psi^s_y\theta(y)^{-1}\, |v^u(x)\rangle\langle v^s_\perp(y)|.
\end{split}
\]
Accordingly,
\begin{equation}\label{eq:jactot}
\begin{split}
J\Phi&:=\det\begin{pmatrix}
            \Id &-\Id\\
            \partial_{x}\eta&
           \partial_{y}\eta 
           \end{pmatrix}=\det\begin{pmatrix}
                         \partial_{x}\eta+
                          \partial_{y}\eta 
                         \end{pmatrix}\\
&=J\Psi^u_x\, J\Psi^s_y\,\theta(x)^{-1}\theta(y)^{-1}\langle v^s_\perp(y),\,
v^u(x)\rangle^2\\
&=J\Psi^u_x\, J\Psi^s_y\,\theta(x)^{-1}\theta(y)^{-1}\det\begin{pmatrix}
                                                         v^u(x)&v^s(y)
                                                         \end{pmatrix}^2.  
\end{split}
\end{equation}
Before starting computing we need to collect some facts.
\begin{sublem}
\label{slem:domains}
If $x\in\partial D_{2\rho}$, $\rho\geq r\geq \cs_4\sqrt\ve$, then, for
$C_4$ large enough,
\begin{enumerate}[i)]
\item if $\|x-y\|\leq \ve$, then $y\not\in D_\rho$.
\item $\|x-\eta\|\leq \frac\ve\rho$.
\item if $\zeta\in W^u(x)$ and $\|\zeta-x\|\leq\frac\ve\rho$, then
$\zeta\not\in D_\rho$.
\end{enumerate}
\end{sublem}
\begin{proof}
The first inequality follows since $D_{2\rho}$ has a vertical size
$4\rho^2$. Thus $2\rho^2-\ve\geq 2\rho^2-\frac
1{\cs_4^2}\rho^2\geq\rho^2$. Such an estimate and Proposition
\ref{prop:dista} imply that the angle between $W^u(x)$ and $W^u(y)$ is
at least $4K_-\rho^{-1}$, thus (ii). Finally, Proposition
\ref{prop:dista} implies that, if $\zeta=(z_1,z_2)$,
$z_2>2\rho^2-\frac{\ve^2}{\rho^2}\geq\rho^2$. 
\end{proof}
We can now start computing the integral.
\[
\begin{split}
\int_{\To^4}dxdy f(x)g(y)q_\ve(x-y)&=\int_{D_{\cs_4 r}^c}dx\int_{\To^2}dy
f(x)g(y)q_\ve(x-y) +\Or(\|f\|_\infty\|g\|_\infty r^3)\\
&\geq\int_{\Phi(D_{\cs_4 r}^c\times \To^2)}\!\!\!\!\!\! d\eta d\xi
f(x)g(y)q_\ve(\xi)J\Phi^{-1}+ 
\|f\|_\infty\|g\|_\infty \Or(r^3)
\end{split}
\]
Next, from formula \eqref{eq:jactot} and Sub-lemma \ref{slem:domains}-$(i)$
follows 
\begin{equation}\label{eq:jac-bound}
|J\Phi(x,y)-1|\leq \cs\frac\ve{\rho^2}.
\end{equation}
Hence,
\[
\begin{split}
\int_{\To^4}dxdy f(x)g(y)q_\ve(x-y)
&=\int_{\To^2} d\eta
f(\eta)g(\eta)+\|f\|_\infty\|g\|_\infty\Or\left(\int_{D_r^c}|1-J\Phi^{-1}|\right)\\  
&\quad+\Or\left(\int_{\Phi(D_{\cs_4 r}^c\times \To^2)} d\eta d\xi
[f(x)-f(\eta)]g(y)q_\ve(\xi)\right) \\
&\quad+\Or\left(\int_{\Phi(D_{\cs_4 r}^c\times \To^2)} d\eta d\xi
f(\eta)[g(y)-g(\eta)]q_\ve(\xi)\right)\\
&\quad +\|f\|_\infty\|g\|_\infty \Or(r^3).
\end{split}
\]
To conclude we must compute the various error terms. 
For each $f\in L^\infty$, holds
\[
\int_{D_R}dx\,f(x)=\frac 23\int_0^Rd\rho\, \frac \rho{1+\rho}\int_{\partial
D_\rho}ds\, f(s).
\]
Remembering
\eqref{eq:jac-bound} and applying Fubini
\[
\int_{D_r^c}|1-J\Phi^{-1}|\leq \cs\|f\|_\infty\int_r^1 d\rho\,
\rho^2\frac\ve{\rho^2}\leq \cs\|f\|_\infty\ve.
\]
Next, let $\gamma_u^\eta:[-\delta,\delta]\to\To^2$ be the unstable manifold
of $\eta$, parametrized by arc-length, and let $s(\eta,\xi)$ be such
that $\gamma_u^\eta(s(\eta,\xi))=x$. Recalling Sub-lemma \ref{slem:domains},
\[
\begin{split}
\left|\int_{\Phi(D_{\cs_4 r}^c\times \To^2)}\!\!\!\!\!\!\!\!\!\!\!\! d\eta d\xi\,
[f(x)-f(\eta)]g(y)q_\ve(\xi)\right|&\leq
\|g\|_\infty\int_{\Phi(D_{\cs_4 r}^c\times \To^2)} \!\!\!\!\!\!\!\!\!\!\!\!d\eta d\xi
\int_0^{s(\eta,\xi)}dt\,|\partial^u f(\gamma_u^\eta(t))|q_\ve(\xi)\\
&\leq\cs \|g\|_\infty\int_{D_{r}^c\times
\To^2} d\eta' d\xi\, |\partial^u f(\eta')|\frac{\|\xi\|}{|\theta(\eta')|}q_\ve(\xi)\\
&\leq \cs \|g\|_\infty\int_0^1 d\rho \;\ve\int_{\partial
D_\rho}|\partial^uf|\\
&\leq  \cs \|g\|_\infty\|\partial^uf\|_{L^1(\nu)}\ve.
\end{split}
\]
Analogously,
\[
\left|\int_{D_{r}^c\times \To^2} d\eta d\xi
f(\eta)[g(y)-g(\eta)]q_\ve(\xi)\right|\leq  \cs
\|f\|_\infty\|\partial^sg\|_{L^1(\nu)}\ve.
\]
We can finally collect all the above estimates and obtain
\[
\begin{split}
\int_{\To^4}dxdy &f(x)g(y)q_\ve(x-y)=\int_{\To^2} d\eta
f(\eta)g(\eta)\\
&+(\|f\|_\infty+\|\partial^uf\|_{L^1(\nu)})(\|g\|_\infty
+\|\partial^sg\|_{L^1(\nu)}) \Or(r^3+\ve)
\end{split}
\]
from which the lemma follows by choosing $r=\ve^{\frac 13}$.
\end{proof}

\section{Decay of correlations}
\label{sec:decay}
Here we put together the results of the previous section to prove
Theorem \ref{thm:main}. 

Let $f,g\in\Co^{1}(\To^2,\R)$, $\int f=0$, then
\[
\begin{split}
\int\Lp^{kn_\ve}f g&=\sum_{i=0}^{k-1}\int\Lp^{n_\ve
i}(\Lp^{n_\ve}-\Lp_\ve)\Lp_\ve^{k-i-1}f g+\int\Lp^k_\ve fg\\
&=\sum_{i=0}^{k-1}\int(\Lp^{n_\ve}-\Lp_\ve)\Lp_\ve^{k-i-1}f g\circ T^{n_\ve
i}+\Or(e^{-\cs k}\|f\|_1\,\|g\|_\infty)\\
&=\sum_{i=0}^{k-1}\int(\Id-Q_\ve)\Lp^{n_\ve}\Lp_\ve^{k-i-1}f g\circ T^{n_\ve
i}+\Or(e^{-\cs k}\|f\|_1\,\|g\|_\infty) ,
\end{split}
\]
where we have used Lemma \ref{lem:l1}. To conclude, by using Lemma
\ref{lem:appr}, we need to estimate the $L^1$ norm of
$\partial^u(\Lp^{n_\ve}\Lp_\ve^{j}f)$. Since $|D\Lp_\ve^{j}f|\leq \cs
\ve^{-1}|f|_\infty$, for $j>0$,
\[
\begin{split}
\|\partial^u(\Lp^{n_\ve}\Lp_\ve^{j}f)\|_{L^1(\nu)}&\leq
\int_0^{\cs_4\sqrt \ve}d\rho\int_{\partial
D_\rho}|\partial^u(\Lp^{n_\ve}\Lp_\ve^{j}f)|
+\int_{\cs_4\sqrt \ve}^1d\rho\int_{\partial
D_\rho}|\partial^u(\Lp^{n_\ve}\Lp_\ve^{j}f)|\\
&\leq \cs\|f\|_\infty \int_0^{\cs_4\sqrt \ve}d\rho\,\ve^{-1}\rho+
\cs|f|_\infty\int_{\cs_4\sqrt \ve}^1d\rho\int_{\partial
D_\rho}\frac{\ve}{\rho^2}\ve^{-1}\\
&\leq \cs \|f\|_\infty\ln\ve^{-1},
\end{split}
\]
where we have used Lemma \ref{lem:expansion} and Sub-lemma
\ref{slem:domains}.\footnote{The above estimate is not sharp. With
some extra work one could avoid the $\ln\ve^{-1}$, yet this would not
change in any substantial way the result, so we chose to keep the
presentation as short as possible.}

Thus
\begin{equation}
\label{eq:sharp}
\begin{split}
\left|\int\Lp^{n}f
g\right|=(\|f\|_{\infty}+\|\partial^u f\|_{L^1(\nu)})\|g\|_{\Co^1}\Or(n\ve^{\frac
32}\ln\ve^{-1}+ 
e^{-\cs n\sqrt \ve}).
\end{split}
\end{equation}
Clearly the best choice is $\ve=\cs (n^{-1}\ln n)^2$ which implies the
Theorem.

\section{Lower bound}
\label{sec:lower}
In this section we prove a lower bound.
\begin{lem}\label{lem:lower}
If there exists a sequence $\gamma_n$ such that, for each
$f,g\in\Co^1$ holds
\[
\left|\int_{\To^2}f\circ T^{-n} g-\int_{\To^2}f\int_{\To^2}g\right| \leq
(\|\partial^u f\|_{L^1(\nu)}+\|f\|_{L^1(\nu)})|g|_{\Co^1}\gamma_n,
\]
then there exists $\cs>0$ such that
\[
\gamma_n\geq n^{-2}\cs.
\]
\end{lem}
\begin{proof}
Let $g\geq 0$ be a smooth function supported away from zero (let us say that
the support of $g$ does not intersect $D_{\frac 12}$). Next, let
$\xi_0=(x_0,y_0)=(\frac 12,y_0)\in W^u(0)$, $\xi_n=T^{-n}\xi_0$. For
each point $\eta$ in a neighborhood of $\xi_n$ let $z_1$ be the
distance, along $W^u(0)$, between $\xi_n$ and $W^u(0)\cap W^s(\eta)$,
and $z_2$ the distance, along $W^s(\xi_n)$, between $\xi_n$ and
$W^u(\xi_n)\cap W^s(\eta)$. By construction 
$\Xi_n(\eta):=(z_1,z_2)$ is a map from a neighborhood of $\xi_n$ to a
neighborhood of the origin with the property that the map transforms
the stable and unstable foliation into the standard foliation given by
the Cartesian coordinates. Clearly, $\Xi_n^{-1}(s,0)=\gamma^u(s)$ (the
unstable manifold of the origin parametrized by arc length and such that
$\gamma^u(0)=\xi_n$), while $\Xi_n^{-1}(0,s)=\gamma_n^s(s)$ (the unstable
manifold of $\xi_n$ parametrized by arc length). Finally we define
$f_n:=(\alpha_n\beta_n)\circ \Xi_n$ with
$\alpha_n(z_1):=\varsigma(\cs_5 n z_1)$, $\beta_n(z_2):=\varsigma(\cs_5n^{-1}
z_2)$, for $C_5$ small enough, and 
\[
\varsigma(x):=\begin{cases}
             1-|x+1|\quad &|x+1|\leq 1\\
             0&|x+1|> 1.
             \end{cases}
\]
In other words, $f_n$ is a function essentially supported on a
neighborhood left of $\xi_n$ of order $n^{-1}$ in the unstable direction
and the stable.  Accordingly, the supports of $f_n\circ T^{-k}$ and $g$
are disjoint for all $k\leq \cs n$. Lemma \ref{lem:holo} implies that
the the support is essentially a rhombus of size $n^{-1}$ and angle
$n^{-1}$.

Thus we have
\[
|g|_{\Co^1}(\|\partial^u f_n\|_{L^1(\nu)}+\|f_n\|_{L^1(\nu)})\gamma_n\geq
\int_{\To^2}f_n\int_{\To^2}g\geq \cs n^{-3}.
\]
On the other hand, using again  Lemma \ref{lem:holo},
\[
\|\partial^u f_n\|_{L^1(\nu)}\leq \cs n\|(\alpha'\beta)\circ
\Xi_n\|_{L^1(\nu)}\leq \cs n^{-1}
\]
which yields the Lemma.
\end{proof}

\begin{rem}
\label{rem:problems} Note that the norms in Lemma \ref{lem:lower} and
in Theorem \ref{thm:main} (even in the stronger version given by
\eqref{eq:sharp}) are different. It is not obvious that, putting the
$L^\infty$ norm instead of the $L^1(\nu)$ one keeps the same rate of
mixing. More generally, it is well known that in the uniformly
hyperbolic setting the smoothness of the function can have an influence
on the mixing rate. An analogous effect may arise in the present
setting but it remains to be investigate. A related problem that needs to be
addressed is the higher dimensional analogous of the present model
where the fixed point has different possibility of losing full
hyperbolicty. It is clear that the present
result is only the starting point and not the end of the story.
\end{rem}

\end{document}